\documentclass[smallextended,onecolumn,envcountsame,envcountsect]{svjour3}       % onecolumn (second format)
\smartqed  % flush right qed marks, e.g. at end of proof
\usepackage[utf8]{inputenc}

\usepackage{graphicx}
\usepackage{mathtools}

\usepackage{amsmath}
\usepackage{amssymb}
\usepackage{amsfonts}
\usepackage{graphicx}
\usepackage{bm}
\usepackage{xcolor}
\usepackage{colortbl}
\usepackage{fixltx2e}
\usepackage{url}
\usepackage{multirow}

\usepackage[colorinlistoftodos,textsize=small]{todonotes}
\newtheorem{observation}[theorem]{Observation}
\usepackage{enumerate}% http://ctan.org/pkg/enumerate

\DeclareMathOperator{\LCA}{LCA}

\newcommand{\leo}[1]{#1}
\newcommand{\celine}[1]{#1}
\newcommand{\todoCeline}[1]{}
\newcommand{\remie}[1]{#1}
\newcommand{\todoRemie}[1]{}
\newcommand{\markj}[1]{#1}

\usepackage{natbib}

\title{Exploring the tiers of rooted phylogenetic network space using tail moves}
\author{Remie Janssen \and Mark Jones \and  P\'eter L. Erd\H{o}s \and  Leo van Iersel \and Celine Scornavacca}

\institute{
             Remie Janssen  \at        
              Delft Institute of Applied Mathematics, Delft University of Technology\\
              Postbus 5031,2600 GA Delft, The Netherlands\\
              \email{remiejanssen@gmail.com}
              \and  
              Mark Jones  \at        
              Delft Institute of Applied Mathematics, Delft University of Technology\\
              Postbus 5031,2600 GA Delft, The Netherlands\\
              \email{markelliotlloyd@gmail.com}
              \and
              P\'eter Erd\H{o}s  \at     
              MTA R{\'e}nyi Institute of Mathematics,\\
			  Re\'altanoda u 13-15\\
			 Budapest, 1053 Hungary\\
              \email{erdos.peter@renyi.mta.hu}
              \and  
               \emph{Corresponding author:}  L. van Iersel \at        
              Delft Institute of Applied Mathematics, Delft University of Technology\\
              Postbus 5031,2600 GA Delft, The Netherlands\\
              \email{l.j.j.v.iersel@gmail.com}
              \and
            C. Scornavacca \at
              Institut des Sciences de l'Evolution, Universit\'{e} de Montpellier, CNRS, IRD, EPHE\\
              Institut de Biologie Computationnelle (IBC)\\
              Place Eug\`ene Bataillon, Montpellier, France\\
              \email{celine.scornavacca@umontpellier.fr}   
}

\begin{document}

\maketitle
%\todo{Indicate what graphics program was used to create the artwork.}
\begin{abstract}
Popular methods for exploring the space of rooted phylogenetic trees use rearrangement moves such as rNNI (rooted Nearest Neighbour Interchange) and rSPR (rooted Subtree Prune and Regraft). Recently, these moves were generalized to rooted phylogenetic networks, which are a more suitable  representation of reticulate evolutionary histories, and it was shown that any two rooted phylogenetic networks of the same complexity are connected by a sequence of either rSPR or rNNI moves. 
Here, we show that this is possible using only \emph{tail moves}, which are a restricted version of  rSPR moves on networks that are more closely related to rSPR moves on trees.
The connectedness still holds even when we restrict to distance-1 tail moves (a localized version of tail-moves).
Moreover, we give bounds on the number of (distance-1) tail moves necessary to turn one network into another, which in turn yield new bounds for rSPR,  rNNI and \remie{SPR} (i.e. the equivalent of \remie{rSPR} on unrooted networks).
%Here, we show that this is also possible using only \emph{tail moves}, which are more restrictive than rSPR moves and in fact more similar to rSPR moves on trees. Moreover, we show that we can even restrict to distance-1 tail moves and we give bounds on the number of (distance-1) tail moves necessary to turn one network into another, which also yield bounds for rSPR and rNNI.
The upper bounds are constructive, meaning that we can actually find a sequence with at most this length for any pair of networks. Finally, we show that finding a shortest sequence \leo{of tail or rSPR moves} is NP-hard.
\end{abstract}

\section{Introduction}

% Rooted phylogenetic networks are a generalisation of  rooted phylogenetic trees where nodes with indegree greater than one are allowed. The difference between the two is that in phylogenetic networks branches may recombine.

Leaf-labelled trees are routinely used in phylogenetics to depict the relatedness between entities such as species and genes.
Accurate knowledge of these trees, commonly-called  \emph{phylogenetic trees}, is vital for our understanding of the processes underlying molecular
evolution, and thousands of phylogenetic trees are reconstructed from molecular data each day. 

However, this representation is not suitable when reticulated events such as hybrid speciations \citep[e.g.][]{Abbott}, horizontal gene transfers \citep[e.g.][]{Zhaxybayeva} and 
recombinations \citep[e.g.][]{Vuilleumier} are involved 
in the evolution of the entities of interest.
In such cases, a more suitable representation can be found in \emph{phylogenetic networks}, where in its broadest sense
a phylogenetic network can be thought of as a leaf-labelled graph (directed or undirected), \leo{usually} without parallel edges and degree-2 nodes 
\citep{Morrison,Huson}. 

Common procedures used to reconstruct phylogenetic trees from biological data are tree rearrangement heuristics \citep{felsenstein2004inferring}. These techniques consist of choosing an optimization criterion (e.g. maximum parsimony,   maximum likelihood, a distance-based scoring scheme, etc.) or opting for a Bayesian approach, and then using \emph{tree rearrangement moves} to explore the space of phylogenetic trees. 
These moves specify possible ways of generating alternative phylogenies
from a given one, and their fundamental property is to be able to transform, by
repeated application, any phylogenetic tree into any other phylogenetic tree. Several tree rearrangement moves have been defined in the past, the most commonly-used ones being  
\emph{NNI} (\emph{Nearest Neighbour Interchange}) moves and  \emph{SPR} (\emph{Subtree Prune and Regraft}) moves; when the phylogenetic trees are considered as \emph{rooted}, i.e. directed and out-branching (i.e. singly rooted), we have their rooted versions: rNNI and rSPR moves. 

\begin{figure}[h!]
\begin{center}
\includegraphics[scale=0.55]{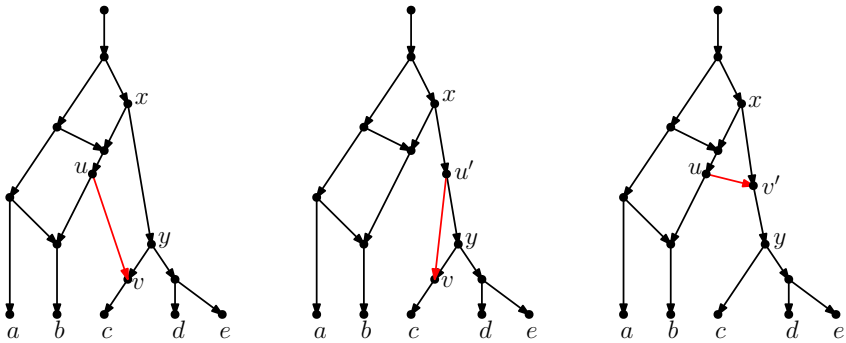}
\end{center}
\caption{
\markj{Starting with the rooted phylogenetic network on the left, if we move the tail of $(u,v)$ to $(x,y)$ we produce the network in the middle. If we instead move the head of $(u,v)$ to $(x,y)$, we produce the network on the right.}
}
\label{fig:NetworkExample}
\end{figure}

Recently, researchers have become interested in defining rearrangement moves for phylogenetic networks and studying their properties
\citep{Bordewich-Lost,francis2017bounds,Gambette-vI-Rearrangement,huber2016transforming}. 
\cite{huber2016transforming} gave a generalization of NNI moves for unrooted phylogenetic networks, and showed the connectivity under these moves of the \emph{tiers} of phylogenetic-network space, i.e. phylogenetic networks having the same \emph{reticulation number}. The latter concept will be formally defined in the next section but, roughly speaking, it is a way to express the amount of reticulate evolution present in a phylogenetic network.  \cite{francis2017bounds} generalized SPR moves for unrooted phylogenetic networks and studied the properties of the NNI and SPR neighbourhoods of a network, giving bounds on their sizes. \cite{Gambette-vI-Rearrangement} focused on \emph{rooted} phylogenetic networks, i.e. phylogenetic networks where the underlying graph is a rooted 
% Leo: I replaced "out-branching" by "rooted" because after some googling I'm still not convinced that out-branching is the right term when talkign about networks
directed acyclic graph, and introduced generalizations of rNNI and rSPR moves for \leo{rooted} phylogenetic networks. rSPR moves consist of \emph{head moves} and \emph{tail moves} (see Figure~\ref{fig:NetworkExample}), while rNNI moves can be seen as distance-1 head moves and distance-1 tail moves \markj{(see Definitions~\ref{def:HeadMove} and~\ref{def:TailMove} for formal \leo{descriptions})}. In the same paper, \citeauthor{Gambette-vI-Rearrangement} showed the connectivity of each tier of rooted phylogenetic networks under rNNI moves (and consequently also under rSPR moves) and gave bounds for the rNNI neighbourhood of a network. Finally, \cite{Bordewich-Lost} introduced another generalisation of rSPR called SNPR (SubNet Prune and Regraft), giving connectivity proofs and bounds for some classes of rooted phylogenetic networks, allowing parallel edges in most cases. %In particular, they proved that each tier of the space of  \emph{tree-child} phylogenetic networks, i.e.  rooted networks in which each internal node has at least one child that does not model a reticulate event, is connected using tail moves only.

Note also that a handful of rearrangement heuristics for phylogenetic networks have \leo{been published} recently (see, for example, PhyloNET \citep{Than} and the BEAST 2 add-on \url{SpeciesNetwork}) and each of them uses its own set of rearrangement moves, e.g. the ``Branch relocator'' move in \citep{Zhang124982} and the ``Relocating the source of an edge'' and ``Relocating the destination of a reticulation edge''  moves in \citep{yu2014maximum}. These sets of moves are often included among the ones cited in the previous paragraph; for example, the  above-cited moves from \citep{yu2014maximum} correspond respectively to   
head and tail moves (see Definitions \ref{def:HeadMove} and \ref{def:TailMove}) and their union correspond to the rSPR \leo{moves} defined in \citep{Gambette-vI-Rearrangement}, c.f. next section. Since these papers do not focus on studying the properties of the moves they define, we will not discuss them here. 

In this paper, we \leo{mostly} focus on rooted phylogenetic networks and \emph{tail moves}. In some sense, these can be seen as the most natural generalisation of rSPR moves on rooted phylogenetic trees to networks. We show that each tier of rooted phylogenetic network space is connected using only tail moves, and even using only distance-1 tail moves. Hence, to get connectivity, head moves are not necessary. We note, however, that head moves could be useful in practice to escape from local optima (see the Discussion section). \leo{We also} analyse the tail-move diameter, giving upper bounds on the number of tail moves, and the number of distance-1 tail moves \leo{(Tail$_1$ moves)}, necessary to turn any rooted phylogenetic network with~$k$ reticulations into any other such network \leo{on the same leaf set}.
% We obtain bounds for the length of sequences of rNNI and rSPR moves that are tighter than the ones given in \citep{Gambette-vI-Rearrangement}.
Since the upper-bound proofs are constructive, we can actually find a sequence to go from one network to another network via tail moves. \leo{Interestingly,} these bounds yield new bounds for rSPR, rNNI and \remie{SPR} moves (\leo{see} Table~\ref{tab:diameters}). %\todo{Leo: I think it would be nice to include  a small table with the bounds.} 

\begin{table}[h]
\centering
\begin{tabular}{|l|c|c|}
\hline
\multicolumn{1}{|c|}{\multirow{2}{*}{\textbf{Move}}} & \multicolumn{2}{c|}{\textbf{Diameter of tier $k$ networks with $n$ leaves}} \\ \cline{2-3} 
\multicolumn{1}{|c|}{}                               & \textbf{Lower bound}                 & \textbf{Upper bound}                 \\ \hline
												     & 										&									   \\
rSPR                                                 & $n-\Theta(\sqrt{n})$                 & $2n+3k$                              \\[1ex] 
Tail                                      			 & $n-\Theta(\sqrt{n})$                 & $3(n+2k)$                            \\[1ex] 
rNNI                                                 & $\Omega(n \log (n))$                 & $O(n^2)$                             \\[1ex] 
Tail$_1$                                             & $\Omega(n \log (n))$                 & $O(n^2)$                             \\[1ex] 
SPR                                                  &                                      & $n+\frac{8}{3}k$                    \\[1ex] \hline
\end{tabular}
\caption{New diameter bounds for several rearrangement moves.}
\label{tab:diameters}
\end{table}
Finally, we show that the computation of a tail move sequence \leo{or rSPR} sequence of shortest length is NP-hard.

%Refer to paper where they use tier for the first time?

% %
% \begin{theorem}[Theorem~3.2 Gambette et al.]\label{the:Gambette}
% Let $N$ and $N'$ be two rooted binary phylogenetic networks on tier $k$ of $X$. Then there exists a sequence of rNNI moves turning $N$ into $N'$. 
% \end{theorem}}

\section{Definitions and properties}
\subsection{Phylogenetic networks}
In this subsection we define the combinatorial objects of interest in this article.

\begin{definition}\label{def:networkTree}
A \emph{rooted binary phylogenetic network} $N=(V,E)$ on a finite set $X$ of taxa is \leo{a directed acyclic graph}
% Leo: I removed out-branching also here
(DAG) with no parallel edges where the \emph{leaves} (nodes of indegree-1 and outdegree-0) are bijectively labelled by~$X$, there is a unique node of indegree-0 and outdegree-1 -- the \emph{root} -- and all other nodes are either  \emph{tree nodes} (indegree-1 and outdegree-2) or  \emph{reticulations} (indegree-2 and outdegree-1). 
We will write $V(N), L(N)$ to denote the nodes and leaves of $N$, respectively. 
\end{definition}
A \emph{rooted binary phylogenetic tree} is a rooted binary phylogenetic network with no reticulation node\leo{s (}all edges are directed outward from the root\leo{)}. For brevity we henceforth simply use the terms \emph{network} and \emph{tree}. 

\begin{definition}
 Let $N$, $N'$ be two networks with leaves labeled with $X$.
 Then an \emph{isomorphism between $N$ and $N'$} is a bijection $\phi: V(N) \rightarrow V(N')$ such that 
 \begin{itemize}
  \item two nodes $u,v \in V(N)$ are adjacent in $N$ if and only if $\phi(u)$ and $\phi(v)$ are adjacent in $N'$;
  \item for any leaf $u \in L(N)$, $\phi(u)$ is the leaf in $L(N')$ that has the same label as $u$.
 \end{itemize}
We say $N$ and $N'$ are \emph{isomorphic} if there exists an isomorphism between $N$ and $N'$.
\end{definition}

Let $u$ and $v$ be nodes in a \leo{network}, then we say $u$ is \emph{above} $v$ and $v$ is \emph{below} $u$ if $u \leq v$ in the order induced by the directed graph underlying the network. Similarly, if $e=(x,y)$ is a directed edge of the network and $u$ a node, then we say that $e$ is \emph{above} $u$ if $y$ is above $u$ and $e$ is \emph{below} $u$ if $x$ is below $u$. 
Let $e=(u,v)$ be an edge of a phylogenetic network, then we say that $u$ is the \emph{tail} of $e$ and $v$ is the \emph{head} of $e$. In this situation, we also say that $u$ is a \emph{parent} of $v$, or $u$ is \emph{directly above} $v$ and $v$ is a \emph{child} of $u$ or $v$ is \emph{directly below} $u$. Note that in a tree $T$, there is always a unique \emph{lowest common ancestor} (LCA) \leo{per pair} of nodes of $T$. This is not the case for networks, where we can have several different LCAs.

A standard measure of tree-likeness of a network $N$ is the \emph{reticulation number}. Denoted by $r(N)$, it is defined as the number of edges that need to be removed in order to obtain a tree.
%More formally:
\markj{As any tree has exactly one more node than edges, we may \leo{equivalently} define it as:}
\[r(N)=\sum_{v\in V : \delta^{-}(v)>0 }\delta^{-}(v)-1 =|E| - |V|+1,\]
where $\delta^{-}(v)$ denotes the indegree of the node $v$. Note that \leo{for binary networks}, $r(N)$ is equal to the number of reticulation nodes.

Here, we are interested in studying the sets of networks with the same reticulation number:

\begin{definition}
Let $X$ be a finite set of taxa. The \emph{$k$-th tier} \citep{francis2017bounds,huber2016transforming} on $X$ is the set of all networks with label set $X$ \leo{and reticulation number~$k$}. % Leo: we have just introduced this term, so we should use it 
\end{definition}

This interest comes from the fact that comparing optimization scores across networks with different reticulation numbers may be tricky, since a network generally allows a better fit with the data when its reticulation number is higher. For this reason, it is increasingly conventional \citep[e.g. in][]{Gambette-vI-Rearrangement} to make a distinction between ``horizontal'' rearrangement moves,  enabling the exploration of a tier, and ``vertical'' moves, allowing a jump across tiers.

The next observation will be useful in the next sections:
%\textcolor{green}{
\begin{observation}\label{obs:numbers}
As every non-root, non-leaf node \markj{in a binary network} is incident to exactly three edges, we have that $3|V| = 2|E| + 2|X| + 2$. Subtracting $2|V| = 2|E| - 2k + 2$ from each side, we get that $|V| = 2|X| + 2k$, which in turn implies that $|E| = 2|X| + 3k - 1$.
Thus any two \markj{binary} networks in the $k$-th tier on $X$ have the same number of nodes and edges, not just the same number of reticulation nodes.
\end{observation}
%} \todo{This is used in a couple of places, and seems basic enough that it's worth mentioning near the start. }

\subsection{Rearrangement moves\label{sec:rearrangement}}

As mentioned in the introduction, several analogues of rearrangement moves on trees have been defined for phylogenetic networks. %\todoCeline{Explain which is which: Gambette: $rNNI=rSPR_1$, Bordewich = SNPR (tail+reticulation!) with parallel edges!, Huber = unrooted NNI}
Such rearrangement moves typically modify the head, the tail, or both the head and the tail of one edge. These moves are \leo{only allowed} if they produce a valid phylogenetic network.

\begin{figure}[h!]
\begin{center}
\includegraphics[scale=0.5]{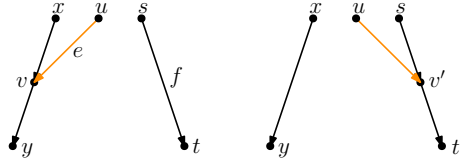}
\end{center}
\caption{A head move as described in Definition~\ref{def:HeadMove}.}
\label{fig:HeadMove}
\end{figure}

\begin{definition}[Head move]\label{def:HeadMove}
Let $e=(u,v)$ and $f$ be edges of a  network. A head move of $e$ to $f$ consists of the following steps:
\begin{enumerate}
\item deleting  $e$;
\item \emph{subdividing} $f$ with a new node $v'$;
\item \emph{suppressing} the indegree-1 outdegree-1 node $v$;
\item adding the edge $(u,v')$.
\end{enumerate}
\emph{Subdividing} an edge $(u,v)$ consists of deleting it and adding a node $x$ and edges $(u,x)$ and $(x,v)$. \emph{Suppressing} an indegree-1, outdegree-1 node $x$ with parent $u$ and child $v$ consists in removing  edges $(u,x)$ and $(x,v)$ and node $x$, and then adding an edge $(u,v)$.

Head moves are only allowed if the resulting digraph is still a network, see Definition \ref{def:networkTree}. 
We say that a head move is a \emph{distance-$d$} move if, after step 2, a shortest path from $v$ to $v'$ in the underlying undirected graph has length $d+1$ (number of edges in the path).
\end{definition}

\begin{figure}[h!]
\begin{center}
\includegraphics[scale=0.5]{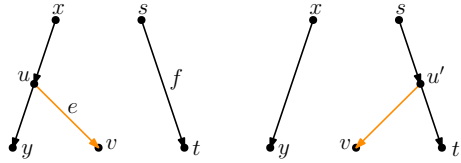}
\end{center}
\caption{A tail move as described in Definition~\ref{def:TailMove}.}
\label{fig:TailMove}
\end{figure}

\begin{definition}[Tail move]\label{def:TailMove}
Let $e=(u,v)$ and $f$ be edges of a  network. A tail move of $e$ to $f$ consists of the following steps:
\begin{enumerate}
\item deleting $e$;
\item subdividing $f$ with a new node $u'$;
\item suppressing the indegree-1 outdegree-1 node $u$;
\item adding the edge $(u',v)$.
\end{enumerate}
Tail moves are only allowed if the resulting digraph is still a network, see Definition \ref{def:networkTree}.
We say that a tail move is a \emph{distance-$d$} move if, after step 2, a shortest path from $u$ to $u'$ in the underlying undirected graph has length $d+1$.
\end{definition}

%Note that connectedness after a move implies that all allowed moves are over a well defined finite distance. Because if this distance were undefined, i.e. no path between $v$ and $v'$ for head moves or $u$ and $u'$ in the case of tail moves, the resulting network would be disconnected.

Note that head moves are only possible for \emph{reticulation edges}, that is, edges in which the head is a reticulation. This means that these moves are, in our opinion, not necessarily part of a natural generalisation of rSPR moves on trees, which consist only of tail moves.

Other generalisations of tree moves that have been proposed  include head moves. For example, one \emph{rSPR move} \citep{Gambette-vI-Rearrangement} on a network consists of one head move or one tail move, and one \emph{rNNI move} consists of one distance-1 head move or one  distance-1 tail move. Thus, it is clear that any tail move is an rSPR move, and any distance-1 tail move is an rNNI move.
\emph{SNPR moves} \citep{Bordewich-Lost} are a variation on the theme: they are defined on networks where parallel edges are allowed, as a tail move or a deletion/addition of an edge. Because the deletion/addition of an edge is a vertical move, SNPR moves can change the reticulation number. Moreover, even if vertical moves are not permitted, the presence of parallel edges makes this restriction of  SNPR moves still subtly different from the tail moves studied in this paper, see Definition \ref{def:TailMove}. \leo{Nevertheless}, Bordewich et al. do disallow parallel edges and vertical moves when studying \emph{tree-child networks} (i.e. networks in which every internal node has at least a child that is not a reticulation) and they prove that any tier of tree-child networks is connected by tail moves if $k\leq |X|-2$. 

\begin{figure}[h!]
\begin{center}
\includegraphics[scale=0.6]{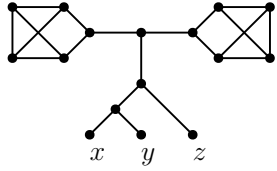}
\end{center}
\caption{\leo{An unrootable unrooted network}.}
\label{fig:HuisjeBoompjeBeestje}
\end{figure}

Rearrangement moves are also defined for \emph{unrooted networks}: connected graphs with \leo{nodes} of degree 1 (the leaves) and of degree 3, where \leo{the} leaves are labelled bijectively with some set $X$\leo{, with~$|X|\geq 2$}. The unrooted equivalents of rSPR moves and rNNI moves are called \emph{SPR moves} and \emph{NNI moves} respectively \citep{huber2016transforming}. An SPR move \leo{relocates} one of the endpoints of an edge, like an rSPR move, with the condition that the resulting graph is still an unrooted network. An NNI move on an unrooted network is again the distance-1 version of an SPR move. \leo{Moreover,} any rSPR move induces an SPR move on the underlying undirected graph, and similarly for rNNI and NNI moves. This means that, in some sense, any rSPR move is an SPR move, and any rNNI move is an NNI move. The converse is clearly not true: for example, an rSPR move that creates a directed cycle is invalid, but the induced SPR move on the undirected network may be valid. In addition, not every unrooted network is the underlying graph of a rooted network (see for example the network in Figure~\ref{fig:HuisjeBoompjeBeestje}). Such networks are unrootable and will be treated in more detail in subsection~\ref{sec:unrooted}.%\todoRemie{A lot of the explanation here would be covered by the material in the new text, make sure we do not have too much redundancy.}

\subsection{Properties of tail moves} 
From now on, we will mostly focus on tail moves (and hence on rooted networks). Note that not all edges of a network can be moved by a tail move; we introduce here the notion of 
\emph{movability} of edges: 

\begin{definition}\label{def:Movable}
Let $e=(u,v)$ be an edge in a  network. Then $e$ is called \emph{non-movable} \leo{if} $u$ is the root, \leo{if~$u$ is} a reticulation, or if the removal of $e$ followed by \leo{suppressing}~$u$ creates parallel edges. Otherwise, $e$ is called \emph{movable}.
\end{definition}

There is only one situation in which \leo{moving a tail of an edge~$(u,v)$} can result in parallel edges: \leo{when there exists an edge from the parent of $u$ to the child of $u$ other than $v$}. The situation is characterized in the following definition:

\begin{definition}\label{def:triangle}
Let $N$ be a  network and let $x,u,y$ be nodes of $N$. We say $x,u$ and $y$ form a \leo{\emph{triangle}} if there are edges $(x,u), (u,y)$ and $(x,y)$. The edge $(x,y)$ is called the \emph{long edge} and $(u,y)$ is called the \emph{bottom edge} \leo{of the triangle}.

\end{definition}
The interesting case in Definition~\ref{def:triangle} is when $u$ is a tree node:
\begin{observation}\label{rem:Triangle}
Let $x,u$ and $y$ form a triangle in a network, and let $v$ be the other child of $u$. The edge $(u,v)$ is \leo{non-}movable because this move would create parallel edges from $x$ to $y$. 
The edges $(u,y)$ and $(x,y)$ are movable, however, %next part can be the lemma
and if \leo{$(u,y)$ is moved sufficiently far up (i.e., destroying the triangle) then the new edge~$(x,v)$ is also movable.} %Leo: note that if you move $(x,y)$ then there won't be a new edge (x,v)!
%
% Here is the old version:
%
% one of those is moved sufficiently far up, the new edge $(x,v)$ is also movable. By sufficiently far, we mean that the move has to break the triangle, resulting in a non-isomorphic network.
\end{observation}
%REMARK LEO:
% it is sufficient to move it to the root edge. if you do that, the edge (x,v) becomes movable, unless the parent of x is the root.
% It might be a good idea to make this a lemma.

The following observation is a direct consequence of Observation~\ref{rem:Triangle} and of the binary nature of the networks studied in this paper, and will play an important role in the arguments presented in the next section:

\begin{observation}\label{rem:MovableTree}
Let $u$ be a tree node, then at least one of its child edges is movable because at most one of these \leo{has its tail} in a triangle \leo{not containing its head}.
\end{observation}

Note that movability of an edge $e$ does not imply that there exists a valid tail move for $e$: it only ensures that we can remove the tail without creating a clear violation of the definition of a network; it does not ensure that we can reattach it anywhere else. The following observation characterizes valid moves:

\begin{observation}\label{rem:MovableTo}
The tail of an edge $e=(u,v)$ can be moved to another edge $f=(s,t)$ if and only if the following conditions hold:
\begin{enumerate}
\item $e$ is movable;
\item $f$ is not below $v$; 
\item $t\neq v$.
\end{enumerate}
\end{observation}

The first condition assures that the tail can be removed, the second that we do not create cycles, and the third that we do not create parallel edges.
Note that these conditions imply that \emph{moving} a (movable) edge \emph{up} is allowed, i.e. moving an edge to another edge that is above it. \leo{In particular, a tail can always be moved to the root edge. However, note that it is not certain that this results in a non-isomorphic network.}

The following lemma is related to the previous observation, and will be used in the next section to find a tail that can be moved down ``sufficiently far''(this concept will be clearer in due time):

\begin{lemma}\label{lem:MovableAncestorOK}
Let $x,y$ be nodes of a phylogenetic network $N$ such that neither~$x$ nor~$y$ is an LCA of~$x$ and~$y$. Then there exists a movable edge in $N$ that is not both above $x$ and above $y$.
\end{lemma}
\begin{proof}
Consider an arbitrary LCA $u$ of $x$ and $y$. This LCA is a tree node, and both child edges are above either $x$ or $y$, but not both. Because at least one of the child edges of a tree node is movable (by Observation~\ref{rem:MovableTree}), at least one of the child edges of the LCA has the desired properties.\qed
\end{proof}

We conclude this section by citing a result by \citeauthor{Gambette-vI-Rearrangement} that implies connectivity of $k$-th tiers (for any $k\geq 0$) via \emph{rNNI moves}, which are equivalent to the combination of distance-1 head moves and  distance-1 tail moves. This result is fundamental for our proof of connectivity of the tiers via tail moves.

%  \textcolor{green}{TO FINISH Here, we study the connectivity of the tiers of phylogenetic network space using only tail moves. We initially prove our connectivity result using the following theorem by Gambette et al.\citep{Gambette-vI-Rearrangement}, who study the \emph{rNNI move} which is equivalent to either a  distance-1 head move or a  distance-1 tail move.
% %
\begin{theorem}[Theorem~3.2 in \citet{Gambette-vI-Rearrangement}]\label{the:Gambette}
Let $N$ and $N'$ be two rooted binary phylogenetic networks belonging to the $k$-th tier on a fixed leaf set $X$. Then there exists a sequence of rNNI moves turning $N$ into $N'$. 
\end{theorem}

%%%%%%%%%%%%%%%%%%%%%%%%%%%%%%%%%%%%%%%%%%%%%%%%%%%%%%%%%%%%%%%%%%%%%%%%%%%%%%%%%%%
%%%%%%%%%%%%%%%%%%%%%%%%%%   SECTION:REWRITING SEQUENCES   %%%%%%%%%%%%%%%%%%%%%%%%
%%%%%%%%%%%%%%%%%%%%%%%%%%%%%%%%%%%%%%%%%%%%%%%%%%%%%%%%%%%%%%%%%%%%%%%%%%%%%%%%%%%
\section{Head moves rewritten}
\begin{figure}[h!]
\begin{center}
\includegraphics[scale=0.5]{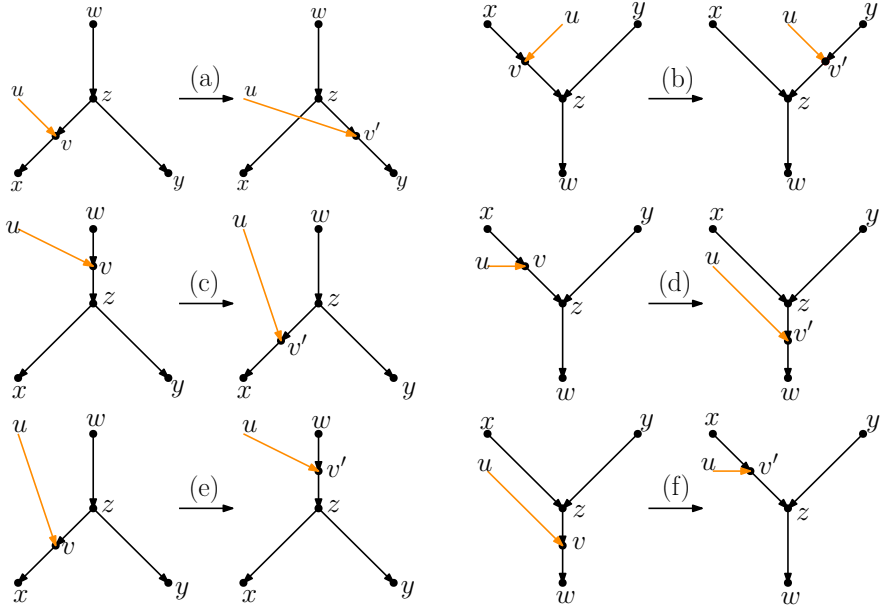}
\end{center}
\caption{Illustration of all possible distance-1 head moves: (a) sideways movement below a tree node, (b) sideways movement above a reticulation, (c) downward movement through a tree node, (d) downward movement through a reticulation, (e) upward movement through a tree node, (f) upward movement through a reticulation.}
\label{fig:HeadMoves}
\end{figure}

In this section we present one of the main results of the paper: the connectivity \leo{of $k$-th} tiers using tail moves. Theorem~\ref{the:Gambette} tells us that $k$-th tiers are connected by rNNI moves. This means that, to prove our result, it suffices to show that any distance-1 head move can be replaced by a sequence of tail moves. To this end, \leo{we list in Figure~\ref{fig:HeadMoves}} all different cases where it is possible to perform a distance-1 head move. We observe the following:

\begin{observation}\label{rem:HeadCases}
Each distance-1 head move is depicted as exactly one of the cases (a)-(f) in Figure~\ref{fig:HeadMoves}. Note that the figure does not indicate whether $u,w,x,y$ are all distinct. There are only a few cases in which $u,w,x,y$ are not all distinct, and the head move is valid. If $x=y$ in cases (a) and (b) or $u=y$ in cases (d) and (f), we have a valid move that results in a network that is isomorphic to the starting one. Having $u=w$ in case (a) is the only situation that results in a non-isomorphic network. All other possibilities lead to invalid moves.
%\todoCeline{I have not checked this in detail but someone should :) }
\end{observation}

\begin{observation}
It is easy to see that moves (c) and (e) as well as moves (d) and (f) are each others reversions. Because all tail moves are also reversible, we only have to show that moves (a)-(d) can be rewritten as a sequence of tail moves.
\end{observation}

We  treat all cases separately in the following lemma. We note here that, in the proof of this lemma, the sequences of tail moves used to mimic the different distance-1 head moves are often non-unique and possibly non-optimal. In the following, we shall use the convention of naming $t'$ the new tail node created when moving an edge with tail $t$.

\begin{lemma}\label{lem:rewriteHeadMoves}
All distance-1 head moves can be substituted by a sequence of tail moves, except for the head move in the network depicted in Figure~\ref{fig:HeadMoveACounter}. For networks with more than one leaf, a sequence of length at most four can be found.
\end{lemma}
%%%%%%%%%%%%%%%%%%%%%%%%%%%%%%%%%%%%%%CASE A
\begin{proof}
We analyse Cases (a)-(d) separately.\\\\
\noindent\fbox{ \parbox{7.5em}{\textbf{Head move $\bm{(a)}$}}}
In light of Observation~\ref{rem:HeadCases}, we distinguish two cases: $u\neq w$ and $u=w$. All other nodes in Figure~\ref{fig:HeadMoves}(a) are distinct.
\begin{enumerate}
\item\label{case:uNotEqualTow} {$\bm{u\neq w}$.} 
In this case we can use the sequence of tail moves depicted in Figure~\ref{fig:HeadMoveA}, since the validity of the head move implies that the intermediate network does not contain directed cycles and parallel edges, as we shall show in the following.

%\remie{
To see that the intermediate network has no directed cycles or parallel edges, we check that the tail move $e$ to $(v,x)$ is valid. Note that $e$ is movable (Definition~\ref{def:Movable}) because $z$ is not the root nor a reticulation, and because $u\neq w$, the removal of the tail of $(z,y)$ does not create parallel edges. Moreover, by  Observation~\ref{rem:MovableTo}, $e$ can be moved to $(v,x)$ since $(v,x)$ is not below $y$ (otherwise there would be a path from $y$ to $v$ in the network after the head move, which implies a directed cycle, and the head move could not be valid) and all nodes in Figure~\ref{fig:HeadMoves} are distinct, so in particular $x\neq y$. 
%\textcolor{red}{Suppose by contradiction that the intermediate network is cyclic; then the directed cycle must involve edge $(z',y)$, and this means that there is a path from $y$ to $z'$, thus a  path from $y$ to $v$. This, in turn, implies that there is a path from $y$ to either $u$ or to $w$. Because the tail move does not change any other part of the network, this path must also exist either in the network before the head move or in the network after the move. This leads to a contradiction, because then one of these networks contains a cycle, making the head move invalid.\\
%It is clear that removing the tail of $(z,y)$ does not create parallel edges because $u\neq w$. The only other possibility for creating parallel edges, is if $x=y$. But, as noted above, in this case the networks before and after the head move are isomorphic.}%%\todoRemie{I added an alternative in green: the green text uses more of the framework that we have built. I thought it would be much shorter, but I disappointed myself.}

Thus, we can conclude that the intermediate network is valid, and therefore that  the sequence of tail moves is valid.

\begin{figure}[h!]
\begin{center}
\includegraphics[scale=0.5]{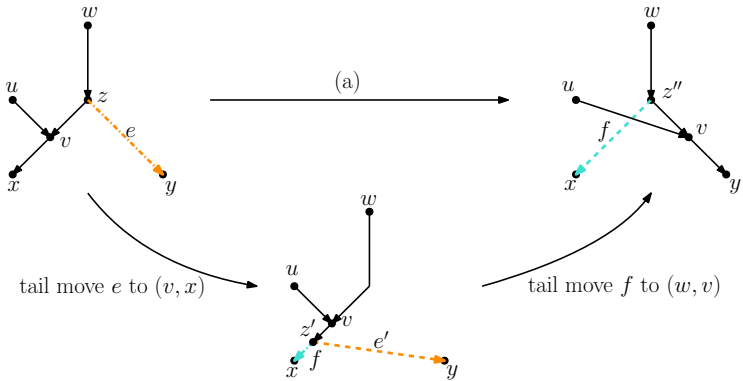}
\end{center}
\caption{\markj{Proof of Lemma~\ref{lem:rewriteHeadMoves}:} The sequence of tail moves needed to simulate head move (a) in Case~\ref{case:uNotEqualTow}. Moving edges are dash-dotted before a move, and dashed after a move.}
\label{fig:HeadMoveA}
\end{figure}

\item\label{case:uEqualsw} {$\bm{u=w}$.} %\todo{MJ - Changed structure of this section (up to head move (b)). Original version is commented out in case we want to switch back.}
Note that $u$ and $v$ form a triangle together with the tree node $z$. In this situation, we cannot directly use the same sequence as before, since this sequence would create parallel edges. There are two conditions under which we can solve this problem:

\begin{enumerate}[i]
\item\label{option:addingTail}{\bf There is a tree node somewhere not above $\bm u$.} 
In this case we ``expand'' the triangle.
Instead of moving edge $e$ directly, we first subdivide it by moving a tail to $e$. Then we can apply the sequence of moves depicted in Figure~\ref{fig:HeadMoveA'}. Barring the addition of the ``extra tail'', the sequence of moves is quite similar to the moves in Case~\ref{case:uNotEqualTow}, and we can prove in a similar way that this will not create cycles or parallel edges.

\begin{figure}[h!]
\begin{center}
\includegraphics[scale=0.5]{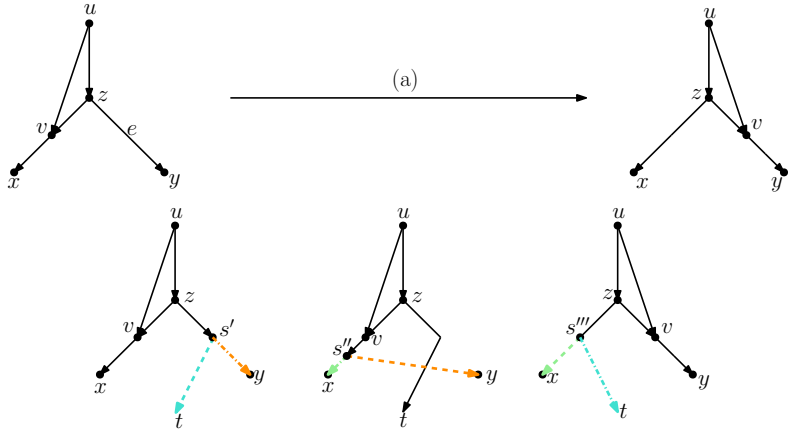}
\end{center}
\caption{\markj{Proof of Lemma~\ref{lem:rewriteHeadMoves}:} The sequence of tail moves used to simulate head move (a) in Case~\ref{case:uEqualsw}. The ``extra tail'' $s$ of edge $(s,t)$ is used in the sequence of moves: $(s,t)$ to $(z,y)$, $(s',y)$ to $(v,x)$, $(s'',x)$ to $(z,t)$ and $(s''',t)$ to $(p,r)$. Here $p$ and $r$ are respectively the parent and (other) child of $s$ in the network before the head move.}
\label{fig:HeadMoveA'}
\end{figure}

\item\label{option:destroyTriangle} {\bf There is at least one vertex above $\bm u$ in addition to the root.}
In this case we ``destroy'' the triangle.
The bottom edge of the triangle $(z,v)$ is movable. If \markj{this edge} is moved to the root \markj{(i.e. to the single edge incident to the root)}, the situation changes to that of Case~\ref{case:uNotEqualTow}. Thus, we can apply the sequence  moves for that case, and move the bottom edge of the triangle back. Since $(z,v)$ is movable, moving up and down to/from the root are valid moves by Observation~\ref{rem:MovableTo}.%\todo{Maybe add a figure?} 
\end{enumerate}

\begin{figure}[h!]
\begin{center}
\includegraphics[scale=0.5]{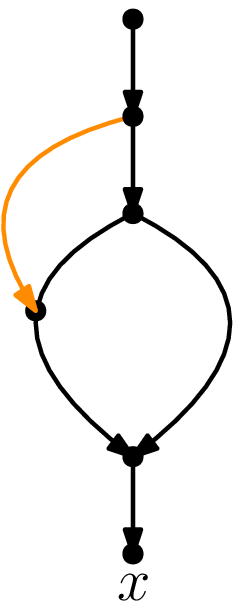}
\qquad\qquad\qquad
\includegraphics[scale=0.5]{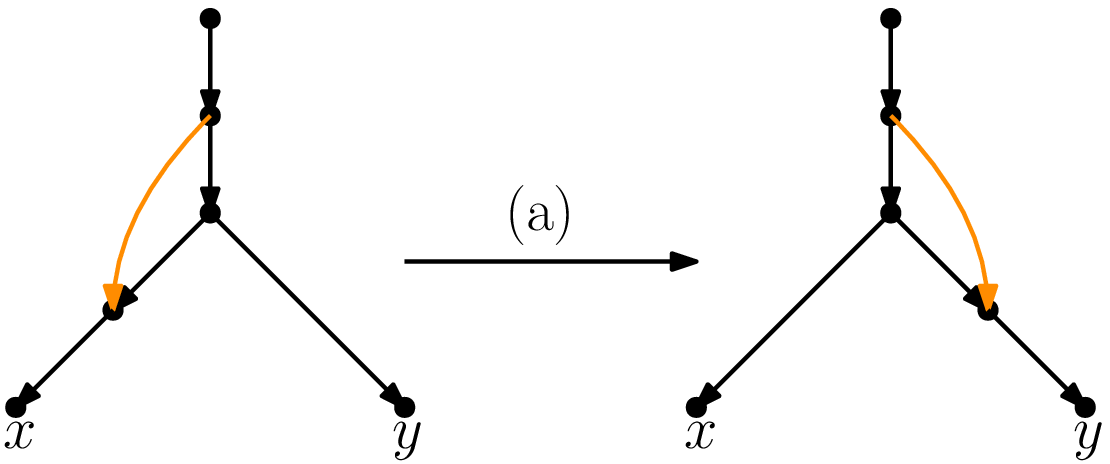}
\end{center}
\caption{The networks for which there are no valid tail moves for simulate a distance-1 head move of type (a). Left: The network with one leaf and 2 reticulations where all valid head moves give isomorphic networks. Right: a head move of type (a) which cannot be substituted by a sequence of tail moves.}
\label{fig:HeadMoveACounter}
\end{figure}

The two networks to which neither of these conditions apply are shown in Figure~\ref{fig:HeadMoveACounter}. The first is the network on one leaf with two reticulations. In this network no head move leads to a different (non-isomorphic) network. The only non-trivial case is the network with two leaves and one reticulation:  the head move cannot be substituted by a sequence of tail moves, because there is no valid tail move.  Note that this network is  excluded in the statement of the lemma.
\end{enumerate}
\begin{figure}[h!]
\begin{center}
\includegraphics[scale=0.5]{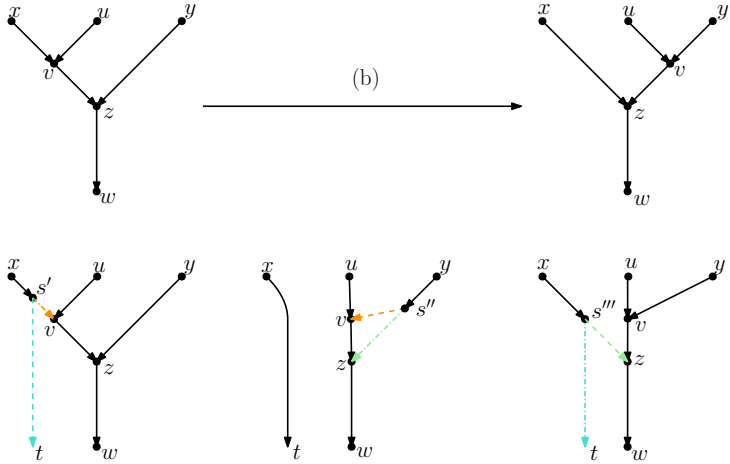}
\end{center}
\caption{\markj{Proof of Lemma~\ref{lem:rewriteHeadMoves}:} The sequence of tail moves used to simulate head move (b). The moving edges are coloured, and the colouring of the edges is consistent throughout the sequence. The tail of the extra edge (light blue) is first moved to $(x,v)$, then the orange edge $(s',v)$ is moved to $(y,z)$, then the green edge $(s'',z)$ is moved to $(x,t)$, and finally the light blue edge $(s''',t)$ (the ``extra tail'') is moved back.}\label{fig:HeadMoveB}
\end{figure}
\noindent\fbox{ \parbox{7.5em}{\textbf{Head move $\bm{(b)}$}}}
The idea of this substitution is to use an extra tail again. The sequence we use is given in Figure~\ref{fig:HeadMoveB}.
The main challenge is to find a usable tail: to use the sequence of moves, we need a tail that can be moved below either $x$ or $y$. We will find this tail by considering an LCA $s$ of $x$ and $y$ in Figure~\ref{fig:HeadMoveB}. We treat the following cases: 
\begin{enumerate}
\item\label{case:BSeparateLCA} {\bf The LCA $\bm s$ is neither equal to $\bm x$ nor to $\bm y$.} At least one of the child edges $(s,t)$ of this LCA can be moved and it is not above both $x$ and $y$ (by Lemma~\ref{lem:MovableAncestorOK}). Hence we can move this tail down below one of these nodes (by Observation~\ref{rem:MovableTo}).
\begin{enumerate}[i]
 \item {\bf The movable edge $\bm{(s,t)}$ is above $\bm y$.} If $t$ is not directly below $x$, we can use the proposed sequence of moves. However, if $t$ is directly below $x$, then $x$ is above $y$, and $x$ is the only LCA of $x$ and $y$. This contradicts our assumption that $s\neq x$ and we can certainly use the sequence of moves depicted in Figure~\ref{fig:HeadMoveB}.
\item\label{case:BSymmetry} 
{\bf The movable edge $\bm{(s,t)}$ is above $\bm x$.} The symmetry of the situation and the reversibility of the sequence of moves lets us reduce to the previous case: We do the sequence of moves in reverse order, where we switch the labels of $x$ and $y$.

\end{enumerate}

% Consider an LCA of $x$ and $y$ in Figure~\ref{fig:HeadMoveB} and suppose that this LCA is different from both $x$ and $y$ (we will see how to deal with this case later). At least one of the child edges of this LCA can be moved and it is not above both $x$ and $y$ (by Lemma~\ref{lem:MovableAncestorOK}). Hence we can move this tail down below one of these nodes (Remark~\ref{rem:MovableTo}).

% Without loss of generality, let the movable edge $(s,t)$ be above $y$ and not above $x$ (if the movable edge is not above $y$, because of the symmetry of the situation and the reversibility of tail moves, we do the sequence of moves in reverse order, where we switch the labels $x$ and $y$). There are two situations in which the sequence of moves is invalid: if $x\neq s$ is a tree node with children $v$ and $t$ (because one of the intermediate networks in the sequence has parallel edges from $x$ to $t$) or if $x=s$ (the sequence cannot be used as shown, because it shows $s'$ and $x$ as distinct nodes). We now give solutions for both these problems.

% \begin{enumerate}
% \item\label{case:parallelxt} {\bf The children of $\bm x$ are $\bm v$ and $\bm t$.}
% Note that we have $(s,t)$ with $s=\LCA(x,y)$ above $y$ movable. Because there is an edge from $x$ to $t$, there is a path from $x$ to $y$. Hence $x$ is the unique LCA of $x$ and $y$, and thus $s=x$. The problematic case `$x\neq s$ is a tree node with children $v$ and $t$' therefore does not occur.

\item\label{case:LCAxyIsx} {$\bm{x}$ \textbf{is an LCA of~$\bm x$ and~$\bm y$}.}
Because $x$ is the LCA of $x$ and $y$ and there is no path from $v$ to $y$ (such a path would imply a path from $z$ to $y$ and thus a cycle in the starting network), $x$ must be a tree node. Let the children of $x$ be $v$ and $t$, and let the parent of $x$ be $p$. Two cases are possible: 
\begin{enumerate}[i]
\item\label{case:LCABothMovable} {\bf $\bm{(x,v)}$ is movable.} In this case we can use a similar sequence as before, but without the addition of a tail: $(x,v)$ to $(y,z)$, $(x',z)$ to $(p,t)$.%\todo{There is no $x'$, it should be $s''$ right? [RJ:If we move the tail of an edge $(\circ,\cdot)$, I have the convention of labelling the new tail node with $\circ'$. I can either define this convention, mention the new label for each move, or keep it implicit.] }
\item\label{case:LCANotMovable} {\bf $\bm{(x,v)}$ is  not movable.} In this case $p,x$ and $t$ form a triangle. We employ a strategy to break the triangle similar to the one used for  Case~\ref{case:uEqualsw} of head move (a).
\begin{itemize}
\item[\textbullet]\label{case:rootTriangle} {\bf There is a node above the triangle besides the root.} 
%\textcolor{blue}{
Moving the long edge of the triangle to the root, we can reduce the problem to Case~\ref{case:LCABothMovable}. Then we apply the same sequence of tail moves and move the long edge back to the original position.%}
\item[\textbullet]{\bf The node above the triangle is the root.}
If there is no node above the triangle, then %\textcolor{blue}{there is a tree node not above $x$: all nodes except $p$ are below $x$ and counting tree nodes and reticulation nodes (the number of tree nodes is $|X|-1+k$), we conclude that there must be at least one tree $r$ node besides $x$ and $p$. Moving one of the child edges $(r,s)$ up to the bottom of the triangle and taking care that the parents of $v$ and $z$ are still distinct, we create the situation of Case~\ref{case:LCABothMovable}. After the corresponding sequence of moves, we move the edge $(r',s)$ back to its original position.}%the parents are only non distinct, i.e. $u=y$ only happens if $par(u)=y$, $r=u$ and we move $(r,s)$ with $s\neq v$, in this case, we can also move $(u,v)$ up, as it does not come from the side of a triangle. 
% \\
%\textcolor{cyan}{
there is a tree node not above $x$, for example an $\LCA$ $s$ of $u$ and $y$. It is easy to see that $s$ cannot be above $x$ since the parent node of $x$ has children $x$ and $t$ and its parent is the root. Note that $s$ can only be a reticulation if it has outdegreee 1, i.e. if $s=u$ or $s=y$. Suppose that $s=u$. Since $v$ and $z$ are reticulations, this would mean that $y$ is below $z$, which implies the existence of a directed cycle in the original network, a clear contradiction. A similar reasoning shows that $s=y$ leads to a contradiction.
Hence, $s$ is a tree node and one of the child edges can be moved up to $(x,t)$. This puts us in the situation of Case~\ref{case:LCABothMovable}, and we can use the corresponding sequence of moves. After the corresponding sequence of moves, we move the edge back to its original position.%}\todo{Blue and Cyan are alternatives here.}
\end{itemize}
\end{enumerate}

\item\label{case:LCAxyIsy} {$\bm{y}$ \textbf{is an LCA of~$\bm x$ and~$\bm y$}.} Like in Case~\ref{case:BSymmetry}, we use the symmetry of the situation and the reversibility of the moves to reduce to the previous case (Case~\ref{case:LCAxyIsx}) by doing the moves in reverse order where the labels of $x$ and $y$ are switched.
\end{enumerate}
In all sequences we gave, no directed cycles were created, because such cycles imply directed cycles in one of the networks before and after the head move. We conclude that any allowed head move of type (b) can be substituted by a sequence of tail moves.

%%%%%%%%%%%%%%%%%%%%%%%%%%%%%%%%%%%%%%CASE C
\noindent\fbox{ \parbox{7.5em}{\textbf{Head move $\bm{(c)}$}}} This is the easiest case, as we can substitute the head move by exactly one tail move: $(z,y)$ to $(w,v)$, with labelling as in Figure~\ref{fig:HeadMoves}. Parallel edges and directed cycles cannot occur because there are no intermediate networks.

%%%%%%%%%%%%%%%%%%%%%%%%%%%%%%%%%%%%%%CASE D

%New approach, always do the tail switch. 
%No need to call it a tail switch
%Cases:
%one of $(u,v)$ and $(y,z)$ is movable: 2 tail moves
%one of $u$ and $y$ is a tree node: harder version: move the triangle with it
%both under reticulation: Subcases
%  2 leaves or more: make tree node using extra tail
%  1 leaf, this leaf is not $w$.
%  no extra tail possible: long way 1 leaf tail switch
\noindent\fbox{ \parbox{7.5em}{\textbf{Head move $\bm{(d)}$}}} This one is the most complicated distance-1 head move to translate in tail moves. We exploit the following symmetry of this case: relabelling $u\leftrightarrow y$ and $v\leftrightarrow z$ transforms the network before the head move into the network after the head move. Let us start with the easiest case:
\begin{enumerate}
\item\label{case:DOneMovable} {\bf Either $\bm{(u,v)}$ or $\bm{(y,z)}$ is movable.}
First assume that $(u,v)$ is movable and that the other child and parent of $u$ are $p$ and $t$, we move $(u,v)$ to $(y,z)$ and then $(u',z)$ to $(p,t)$. The case that $(y,z)$ is movable can be tackled in a similar way thanks to the reversibility of tail moves and the symmetry given above.
\item\label{case:DTriangle} {\bf Either $\bm{u}$ or $\bm{y}$ is a tree node.}
If one of them is movable, we are in the previous situation, otherwise we mimic the approach of Case~\ref{case:LCANotMovable} of head move (b). Assume without loss of generality that $y$ is at the side of a triangle. We either move the bottom edge of the triangle up to the root, or we use a child edge of a tree node to subdivide the bottom edge of the triangle. %\textcolor{blue}{We can always do this if the triangle is directly below the root: there is a tree node below $y$ that we can use for the same reasons as in Case~\ref{case:rootTriangle}}.\textcolor{cyan}{
In the latter case, we can take, for example, $\LCA(x,u)$, which is a tree node because $x$ and $u$ are both directly above $v$, and which is not above $y$.
%}
%\todo{colours correspond to the previous coloured text.}

%\todo{Again, you should prove more formally that there is always a tree node below $y$.}

% \todo{There is an option which does not break the triangle, and may be shorter if we consider tail$_1$ moves. Maybe it is good to place a remark somewhere about these kinds of sequences. Most of the sequences we use are non-unique and possibly also non-optimal.}
% A third option, which does not break the triangle (and is hence shorter, when we convert to tail_1 moves, is given in Figure~\ref{fig:HeadMoveDTriangle}. It only works if $t$ is not above $u$, otherwise we create cycles!
%
% \begin{figure}[h!]
% \begin{center}
% \includegraphics[scale=0.5]{HeadMoveDTriangle}
% \end{center}
% \caption{The sequence of tail moves used to simulate head move (d) in Case~\ref{case:DTriangle}.}
% \label{fig:HeadMoveDTriangle}
% \end{figure}

\item {\bf Both $\bm{u}$ and $\bm{y}$ are reticulations.}
In this case we try to recreate the situation of Case~\ref{case:DOneMovable} by adding a tail on one of the edges $(u,v)$ and $(y,z)$.
\begin{enumerate}[i]
\item\label{case:HeadDTwoLeaves} {\bf The network contains at least two leaves.}
We can pick two distinct leaves, at least one of which is below $w$. Suppose first that both of these leaves are below $w$, then an LCA of the leaves is also below the reticulation and one of its child edges can be used as the extra tail by Lemma~\ref{lem:MovableAncestorOK}. 
If only one of the leaves is below the reticulation, then any LCA of the leaves has one edge that is not above $w$. If this edge can be moved, we can directly use it, otherwise, we first move the lower part of the triangle to the root, and then still use this edge.%Can we move this one to the root? Is it not below the root already? Explain more carefully
%Indeed might not be possible, but then more tree nodes below the triangle, take LCA(u,y) instead?

\item {\bf The network contains one leaf, and $\bm{w}$ is not this leaf.}
The only leaf is below the lowest reticulation $r$ of the network, which is not $z$. %\textcolor{blue}{
Let the parents of $r$ be $p_1$ and $p_2$. To find a usable extra tail, we use the same argument as in Case~\ref{case:HeadDTwoLeaves} but now with $p_1$ and $p_2$ instead of the two leaves. Note that $\LCA(p_1,p_2)$ is a tree node because both $p_1$ and $p_2$ are directly above $r$.%\todoCeline{I removed ``and there is no path from $r$ to neither $p_1$ nor $p_2$.'' because I think it is not needed}%}%\todo{this is the third time I use this argument to show that an LCA is a tree node, maybe we can incorporate this in the Lemma.}

% \textcolor{red}{By removing the lowest reticulation and the leaf below it, and then adding leaves to the loosened edges%\todo{This is weird because it is NOT a tail move so we cannot say that the space is connected via tail moves. [RJ. It is a tail move, I just use a strange argument here: I guess I was too lazy to rewrite the arguments before using the parents of the lowest reticulation instead of the two leaves. Now I realise that is what I do, I could just write it like that, might be less confusing.] %We may want to not consider networks with one leaf, This will also avoid case 3c, which is not elegant. [RJ. Here I agree: it is not elegant. However, it feels more complete to also consider one leaf networks, as it also works for those networks, be it with an exceptional case 3c.]}
% , we reduce to Case~\ref{case:HeadDTwoLeaves}. The associated sequence of moves where we treat the parents of $r$ as leaves substitutes the head move.}\todoCeline{I still do not think that it is a tail move}

\item {\bf The only leaf in the network is $\bm{w}$.}\label{case:DOneLeaf}
We have not found a way to solve this case `locally' as we did before. Go up from these reticulations to some nearest tree node. The idea, then, is to use one tree edge to do the switch as in the case we discussed previously (Figure~\ref{fig:TailSwitchHard}). 

The result is that all reticulations on the path to the nearest tree node move with the tree edge. These reticulations can be moved back to the other side using the previously discussed moves. In particular we move the heads sideways using head move (b), and then we move them up using head move (f) where the main reticulation ($z$) is not the lowest reticulation in the network. \qed
%
%Hard Tail Switch! More details
%If the tree cannot be moved, move the edge above it, then move this tree edge back. The triangle is restored.
%
%
%
%
%

\begin{figure}[h!]
\begin{center}
\includegraphics[scale=0.4]{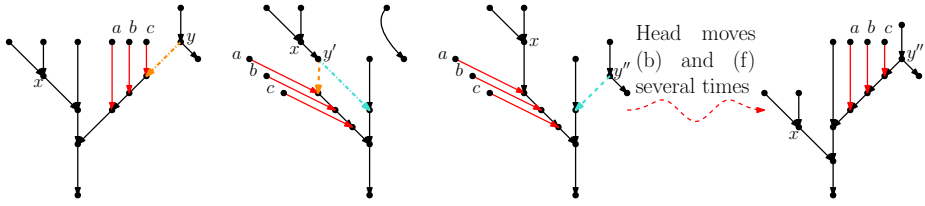}
\end{center}
\caption{\markj{Proof of Lemma~\ref{lem:rewriteHeadMoves}: The sequence of tail moves used to simulate head move (d) in Case~\ref{case:DOneLeaf}.} All reticulations on the path to the nearest tree node above it move with it to the other side.}
\label{fig:TailSwitchHard}
\end{figure}

\end{enumerate}
\end{enumerate}
\end{proof}

Using Theorem~\ref{the:Gambette} and the above lemma, we directly get the connectivity result for tail moves.
\begin{theorem}\label{the:HeadToTail}%Maybe place this earlier, at first mention?
Let $N$ and $N'$ be two networks in the $k$-th tier on~$X$. Then there exists a sequence of tail moves turning $N$ into $N'$, except if $k=1$ and $|X|=2$. 
\end{theorem}

\section{Distances and diameter bounds}
% We use the notation from Remie's notes.
% In particular, we assume that phylogenetic networks are binary, and have a single root of outdegree $1$, indegree $0$.

Given a class of moves $\mathcal{M}$ and two networks $N,N'$, define the \emph{$\mathcal{M}$-distance} $d_{\mathcal{M}}(N,N')$ to be the minimum length of a sequence of moves in $\mathcal{M}$ turning $N$ into $N'$, or $\infty$ if no such sequence exists. 
\celine{Also, denote by $\Delta_k^{\mathcal{M}}(n)$ the \emph{diameter} of the $k$-th tier  with $n$ leaves with respect to $\mathcal{M}$; that is, the maximum value of $d_{\mathcal{M}}(N,N')$ over any pair of networks $N,N'$ belonging to the $k$-th tier on a fixed leaf set $X$ with $|X|=n$. }
In this section, we \leo{mostly} focus on the tail-distance, studying its diameter and the complexity of computing such a distance between two networks.
Interestingly enough, the findings on the tail-distance will also yield results for the rNNI-, rSPR- and SPR-distance.
\celine{Indeed, tail moves are related to other types of moves such as rNNI and rSPR, which has implications for their induced metrics and diameters. There are obvious bounds when a class of moves contains another class of moves, by which we mean that each move of the second class is equivalent to a move of the first class \markj{(e.g. $\text{rNNI}\subseteq\text{rSPR}$)}. Then we have that $\mathcal{M}\subseteq \mathcal{M}'$ implies $\Delta_k^{\mathcal{M}}(n)\geq \Delta_k^{\mathcal{M}'}(n)$. This observation will be very useful in the \markj{rest of this} section.}

\subsection{The complexity of computing the tail distance and the rSPR distance}%\todoCeline{This section contains a lot of non introduced concepts and I do not know how to make it self-contained without making the paper awfully long. Leo: proposed solution is to add a definition of agreement forests and to move the cluster- and subtree reductions to the discussion section, without precise definition.}
In this subsection we prove that computation of the tail and \leo{r}SPR distance between two networks is NP-hard:

\begin{theorem}\label{the:NPhard}
Computing the rSPR distance and the tail distance between two networks is NP-hard for any tier of phylogenetic network space.
\end{theorem}

To prove the theorem, we shall introduce several new concepts.
% \begin{definition}%tree with biconnected component above it
% A \emph{mycorrhizal tree} with leaves $X$ is a network $M=(T,R)$ defined by a tree $T$ with leaves $X$ \textcolor{green}{(the \emph{tree component})} and a biconnected network $R$ with one leaf \textcolor{green}{(the \emph{root component})}. The mycorrhizal tree is the network where the root edge of $T$ is identified with the leaf edge of $R$.
% If $M$ and $M'$ are mycorrhizal trees with trees $T$  and $T'$, then we denote $d_{\text{treeSPR}}(M,M'):=d_{\text{rSPR}}(T,T')$.
% \end{definition}

\begin{figure}[h!]
\begin{center}
\includegraphics[scale=0.5]{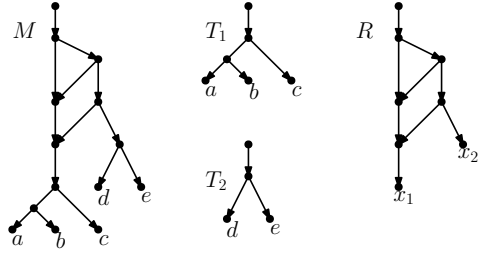}
\end{center}
\caption{From left to right: a mycorrhizal forest $M=(\mathcal{T},R)$, its tree components $\mathcal{T} = (T_1,T_2)$ and the root component $R$.}%the tree components $\mathcal{T} = (T_1,T_2)$ of $M=(\mathcal{T},R)$; the root component $R$ of $M=(\mathcal{T},R)$, the mycorrhizal forest $M=(\mathcal{T},R)$.}
\label{fig:mycorrhzialForest}
\end{figure}

\begin{definition}%tree with biconnected component above it
A \emph{mycorrhizal forest}\footnote{\markj{
 Mycorrhizal forests are so-called because of their similarity to real-life mycorrhizal networks, in which a number of trees may be connected together by an underground network of fungi.}} with leaves $X$ is a network $M=(\mathcal{T},R)$ defined by:
\begin{itemize}
\item a set of $t$ trees $\mathcal{T} = (T_1,\ldots ,T_t)$ with leaf sets~$X_1,\ldots ,X_t$ respectively (the \emph{tree components}),
\item a network $R$ (the \emph{root component}) with $t$ leaves~$x_1,\ldots,x_t$, such that deleting all these leaves and the root makes~$R$ biconnected.
\end{itemize}
The mycorrhizal forest is the network where the root edge of $T_i$ is identified with the edge leading to leaf $x_i$ in $R$ (see Figure~\ref{fig:mycorrhzialForest} for an example). A mycorrhizal forest with~$t=1$ is called a \emph{mycorrhizal tree}.

If $M$ and $M'$ are mycorrhizal forests with tree components~$T_1,\ldots ,T_t$ and~$T'_1,\ldots ,T'_t$, respectively, both having leaf sets $X_1,\ldots ,X_t$, then we denote $d_{\text{treeSPR}}(M,M'):=\sum_{i}d_{\text{rSPR}}(T_i,T_i')$, which is the distance induced by rSPR moves only within tree components (treeSPR moves).
\end{definition}

\begin{definition}
Let $T$ be a tree on $X$ and let $Y\subseteq X\cup\{\rho\}$, where $\rho$ denotes the root of $T$. Then $T|_Y$ is the subtree of $T$ induced by $Y$: Take the union of all shortest paths between nodes of $Y$, and then suppress all indegree-1 outdegree-1 vertices. 
\end{definition}

\begin{definition}
Let $\mathcal{T}=\{T_i\}_{1\leq i\leq t}$ be a set of trees with labels $X$. An \emph{agreement forest} (AF) of $\mathcal{T}$ is a partition $\{X_j\}$ of $X\cup\{\rho\}$ (where $\rho$ denotes the root), such that all $T_i|_{X_j}$ are isomorphic for each fixed $j$ and all $T_i|_{X_j}$ are node-disjoint for a fixed $i$.
\end{definition}

% \begin{lemma}\label{lem:distanceRelationMycorrhizal}
% Let $M$ and $M'$ be two mycorrhizal trees with the same root component and tree \textcolor{green}{components} $T$ and $T'$. Then 
% \[d_{\text{tail}}(M,M')=d_{\text{rSPR}}(M,M')=d_{\text{treeSPR}}(M,M').\]
% \end{lemma}
% \begin{proof}
% Clearly $d_{rSPR}\leq d_{treeSPR}$ and $d_{tail}\leq d_{treeSPR}$, because a treeSPR sequence \todo{define treeSPR sequence} is also a tail sequence and a rSPR sequence.

% Now suppose we have a sequence of tail moves or SPR moves from $M$ to $M'$. Because both networks have all reticulations in the biconnected root component, deleting all moving edges gives a deletion agreement forest on the trees $T$ and $T'$. This deletion agreement forest has size larger than \textcolor{green}{or equal to} the maximal agreement forest (MAF). Because the SPR distance is equal to the size of the MAF (\citep{bordewich2005computational}), we conclude that the number of moves in our sequence is larger than \textcolor{green}{or equal to} the needed number of moves in a treeSPR sequence between $T$ and $T'$. Hence we conclude $d_{rSPR}\geq d_{treeSPR}$ and $d_{tail}\geq d_{treeSPR}$.
% \qed\end{proof}

\begin{lemma}\label{lem:distanceRelationMycorrhizal}
Let $M$ and $M'$ be two mycorrhizal forests with the same root component and tree sets $\{T_i\}_{1\leq i\leq n}$ and $\{T'_i\}_{1\leq i\leq n}$ with leaf sets $X_i=X_i'$. Then 
\[d_{\text{tail}}(M,M')=d_{\text{rSPR}}(M,M')=d_{\text{treeSPR}}(M,M').\]
\end{lemma}
\begin{proof}
%\todo{RJ: I do not think the proof uses that the root component is biconnected: the proof only needs the networks after ``tree cluster reduction'' to be the same. Rewriting will make it a bit harder, because then the number of pendant trees is not in the definition of the network any more, so still prefer the biconnected components.}
Clearly $d_{rSPR}\leq d_{treeSPR}$ and $d_{tail}\leq d_{treeSPR}$, because a treeSPR sequence is also a tail sequence and a rSPR sequence.

Now suppose we have a sequence of tail moves or \leo{r}SPR moves from $M$ to $M'$. Because both networks have all reticulations in the % biconnected
\leo{root component}, deleting all moving edges gives 
% deletion 
agreement forests
%\todoRemie{Is it clear to everyone what a deletion agreement forest is, or should we define this separately? {\bf Celine}: Just remove the word "deletion"and it would be clear.}
on the trees $T_i$ and $T_i'$ for each $i$. These 
% deletion 
agreement forests each have size larger than or equal to the corresponding \leo{maximum} agreement forest (MAF). Because the rSPR distance is equal to the size of \leo{a} MAF \citep{bordewich2005computational} \leo{minus one}, we conclude that the number of moves in the rSPR or tail sequence is larger than or equal to the needed number of moves in all treeSPR sequence\leo{s} between $T_i$ and $T_i'$ together. Hence we conclude $d_{rSPR}\geq d_{treeSPR}$ and $d_{tail}\geq d_{treeSPR}$.
\qed\end{proof}

Theorem~\ref{the:NPhard} follows directly from this lemma, because computing treeSPR distance is NP-hard \citep{bordewich2005computational}. Indeed, for any $k$ we can let $R$ be a network with $k$ reticulations \leo{that becomes biconnected after deleting the root and all leaves}. Then calculating the \leo{r}SPR or tail distance between $M = (\mathcal{T},R)$ and $M' = (\mathcal{T}',R)$ is equivalent to calculating the rSPR distance between $T_i$ and $T_i'$ for each $i$, and $M$ and $M'$ are in the \leo{$k$-th tier}.

\subsection{The diameter of tail \celine{and rSPR} moves}

In this subsection, we  study the diameter bounds of tail move operations, i.e. the maximum value of $d_{Tail}(N,N')$ over all possible networks $N$ and $N'$ with the same reticulation number. 
In the following, we shall use the convention of naming $t_1$ the new tail node created when moving an edge with tail $t$, and naming $t_{i+1}$ the node created when moving an edge with tail $t_i$.
%\todo{Added a note about the notation convention in this section (I thought this was better than rewriting everything to match previous section, which would probably introduce errors) - MJ}

Given a network $N$ and a set of nodes $Y$ in $N$,
we say $Y$ is \emph{downward-closed} if for any $u \in Y$, every child of $u$ is in $Y$.

\begin{lemma}\label{lemma:GreenLineLemma}
 Let $N$ and $N'$ be networks 
 %\textcolor{green}{
 in the $k$-th tier on $X$
 %}
 %on $X$ with the same number of reticulations,
such that $N$ and $N'$ are not the networks depicted in Figure~\ref{fig:HeadMoveACounter}.
 Let $Y \subseteq V(N), Y' \subseteq V(N')$ be downward-closed sets of nodes such that
 $L(N)\subseteq Y, L(N') \subseteq Y'$, and $N[Y]$ is isomorphic to $N[Y']$.
 Then there is a sequence of at most $3|N\setminus Y|$ tail moves turning $N$ into $N'$.
\end{lemma}
\begin{proof}
%\textcolor{red}{By setting $Y = L(N)$ and $Y' = L(N')$, we note that such $Y$ and $Y'$ always exist.}\todo{This sentence is unclear because in the following of the proofs, these equivalences do not necessarily hold.}
We first observe that any isomorphism between $N[Y]$ and $N[Y']$ maps reticulations (tree nodes) of $N$ to reticulations (tree nodes) of $N'$.
Indeed, every node in $N$ is mapped to a node in $N'$ of the same outdegree, and the tree nodes are exactly those with outdegree $2$.
It follows that $Y$ and $Y'$ contain the same number of reticulations and the same number of tree nodes. 
%Furthermore, $N$ and $N'$ have the same number of nodes. Indeed, we have that  $|E(N)| = (3|V(N)| - 2|X| - 2)/2$ and similarly $|E(N')| = (3|V(N')| - 2|X| - 2)/2$ (as every node has combined in- and outdegree $3$, except for the leaves and the root).
%But 
%as $N$ and $N'$ have the same reticulation number, we also have that $|E(N)| - |V(N)| = |E(N')| - |V(N')|$.
%This implies that $(|V(N)| - 2|X| - 2)/2 = (|V(N')| - 2|X| - 2)/2$, and hence $|V(N)| = |V(N')|$.
%\todo{MJ - cut explanation why $N,N'$ have same number of nodes and edges (now in Section 2)}
As $N$ and $N'$ have the same number of nodes and reticulations by Observation~\ref{obs:numbers}, it also follows that $V(N)\setminus Y$ and $V(N')\setminus Y'$ contain the same number of reticulations and the same number of tree nodes.
%as $N$ and $N'$ have the same number of reticulations (\textcolor{red}{and thus the same number of tree nodes}\todo{you need to prove it}), it also follows that $V(N)\setminus Y$ and $V(N')\setminus Y'$ contain the same number of reticulations and the same number of tree nodes.

 We prove the claim by induction on $|N \setminus Y|$.
 If  $|N \setminus Y| = 0$, then $N = N[Y]$, which is isomorphic to $N'[Y'] = N'$, and so there is a sequence of $0$ moves turning $N$ into $N'$.
 
 If  $|N \setminus Y| = 1$, then as $Y$ is downward-closed, $N \setminus Y$ consists of $\rho_N$, the root of $N$, and by a similar argument $N '\setminus Y'$ consists of $\rho_{N'}$.
 Let $x$ be the only child of $\rho_N$, and note that in $N[Y]$, $x$ is the only node of indegree $0$, outdegree $2$.
 It follows that in the isomorphism between $N[Y]$ and $N'[Y']$, $x$ is mapped to the only node in $N'[Y']$ of indegree $0$, outdegree $2$, and this node is necessarily the child of $\rho_{N'}$.
 Thus we can extend the isomorphism between $N[Y]$ and $N'[Y']$ to an isomorphism between $N$ and $N'$ by letting $\rho_N$ be mapped to $\rho_{N'}$.
 Thus again there is a sequence of $0$ moves turning $N$ into $N'$.
 
 So now assume that $|N \setminus Y| > 1$.
 We consider three cases, which split into further subcases.
 In what follows, a \emph{lowest node of $N \setminus Y$ (or $N' \setminus Y'$)} is a node $u$ in $V(N)\setminus Y$ (or $V(N')\setminus Y'$) such that all descendants of $u$ are in $Y$ ($Y'$).
 Note that such a node always exists, as $Y$ ($Y'$) is downward-closed and $L(N)\subseteq Y$ ($L(N') \subseteq Y'$).
 
 \begin{enumerate}
  \item\label{case:reticulation}  {\bf There exists a lowest node $\bm u'$ of $\bm{N' \setminus Y'}$ such that $\bm{u'}$ is a reticulation:} In this case, let $x'$ be the single child of $u'$.
 Then $x'$ is in $Y'$, and therefore there exists a node $x \in Y$ such that $x$ is mapped to $x'$ by the isomorphism between $N[Y]$ and $N'[Y']$.
 Furthermore, $x$ has the same number of parents in $N$ as $x'$ does in $N'$ (because the networks are binary and  $x$ has the same number of children in $N$ as $x'$ does in $N'$ by the isomorphism between $N[Y]$ and $N'[Y']$), and the same number of parents in $Y$ as $x'$ has in $Y'$.
 Thus $x$ has at least one parent $z$ such that $z$ is not in $Y$.
 
 We now split into two subcases:
 
    \begin{enumerate}
    \item\label{case:reticulationFound} {\bf $\bm z$ is a reticulation:} in this case, let $Y_1 = Y \cup \{z\}$ and $Y_1' = Y' \cup \{u'\}$, and extend the isomorphism between $N[Y]$ and $N[Y']$ to an isomorphism between $N[Y_1]$ and $N[Y_1']$, by letting $z$ be mapped to $u'$ (see Figure~\ref{fig:diamReticEasy}).
    % \todo{ATTENTION! In Fig 10 you have that $v$ labels both a parent of $u$ and itis child! CHANGE}
    We now have that $Y_1$ and $Y_1'$ are downward-closed sets of nodes such that $N[Y_1]$ is isomorphic to $N[Y_1']$, and $L(N)\subseteq Y_1, L(N') \subseteq Y_1'$.
    Furthermore $|N \setminus Y_1| = |N \setminus Y|-1$.
    Thus by the inductive hypothesis, there is a sequence of $3|N \setminus Y_1| =  3|N \setminus Y|-3$ tail moves turning $N$ into $N'$.
    
    \begin{figure}[h!]
\begin{center}
\includegraphics[scale=0.5]{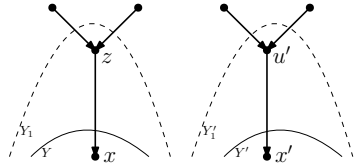}
\end{center}
\caption{
%Illustration of Case~\ref{case:reticulationFound}.
\markj{Proof of Lemma~\ref{lemma:GreenLineLemma}, Case~\ref{case:reticulationFound}:}
If $u'$ is a lowest reticulation in $N' \setminus Y'$ with child $x'$, and the node $x \in Y$ corresponding to $x'$ has a reticulation parent $z$ in $N \setminus Y$, then we may add $z$ to $Y$ and $u'$ to $Y'$. }
\label{fig:diamReticEasy}
\end{figure}
 
    \item {\bf $\bm z$ is not a reticulation:} then $z$ cannot be the root of $N$ (as this would imply  $|N \setminus Y| = 1$), so $z$ is a tree node. It follows that the edge $(z,x)$ is movable, unless the removal of $(z,x)$ followed by suppressing $z$ creates parallel edges.

	\begin{enumerate}
	 \item\label{case:reticulationMoveEasy} {\bf $\bm{(z,x)}$ is movable:} In this case, let $u$ be any reticulation in $N \setminus Y$ (such a node must exist, as $u'$ exists and $N \setminus Y$,  $N' \setminus Y'$ have the same number of reticulations).
	 Let $v$ be the child of $u$ (which may be in $Y$), and observe that the edge $(u,v)$ is not below $x$ (as $x \in Y$ and $u \notin Y$). If $v = x$, then $u$ is a reticulation parent of $x$ that is not in $Y$, and by substituting $v$ for $z$, we have case~\ref{case:reticulationFound}. So we may assume $v \neq x$. Then it follows from Observation~\ref{rem:MovableTo} that the tail of $(z,x)$ can be moved to $(u,v)$.	 
	 Let $N_1$ be the network derived from $N$ by applying this tail move, and let $z_1$ be the new node created by subdividing $(u,v)$ during the tail move (see Figure~\ref{fig:diamReticMoveEasy}). 
     %\todo{Mark, please swap $v$ and $y$ in Fig 10, you use $v$ as child of $u$ in the text!} 
     Thus, $N_1$ contains the edges $(u,z_1), (z_1,v), (z_1,x)$.
     (Note that if $z$ is immediately below $u$ in $N$ i.e. $z=v$, then in fact $N_1 = N$. In this case we may skip the move from $N$ to $N_1$, and in what follows substitute $z$ for $z_1$.)
     %What if $v=z$? I guess it still works: let $y$ be the other child of $z$, move $(z,y)$ up. This is possible because: it cannot create parallel edges ($(u,x)$ is not an edge), and there is something above $u$ (it is a reticulation).
     
     Note that $(z_1,v)$ is movable in $N_1$, since the parent $u$  of $z$ is a reticulation node, and therefore deleting $(z_1,v)$
     %$u$ \todo{why deleting $u$?} 
     and suppressing $z_1$ cannot create parallel edges. Let $w$ %$a$ 
     be one of the parents of $u$ in $N_1$. Then the tail of $(z_1,v)$ can be moved to $(w,u)$ (as $u \neq v$, and $(w,u)$ is not below $v$  as this would imply a cycle in $N_1$).
    
     So now let $N_2$ be the network derived from $N_1$ by applying this tail move (again see Figure~\ref{fig:diamReticMoveEasy}). In $N_2$, the reticulation $u$ is the parent of $x$ (as $z_1$ was suppressed), and thus Case~\ref{case:reticulationFound} applies to $N_2$ and $N'$. Therefore there exists a sequence of $3|N \setminus Y|-3$ tail moves turning $N_2$ into $N'$. As $N_2$ is derived from $N$ by two tail moves, there exists a sequence of $3|N \setminus Y|-1$ tail moves turning $N$ into $N'$.
     
         \begin{figure}[h!]
\begin{center}
\includegraphics[scale=0.5]{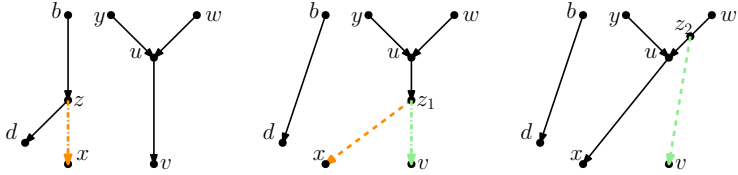}
\end{center}
\caption{
%Illustration of Case~\ref{case:reticulationMoveEasy}. 
\markj{Proof of Lemma~\ref{lemma:GreenLineLemma}, Case~\ref{case:reticulationMoveEasy}:}
If $(z,x)$ is movable, we may move the tail of $(z,x)$ to $(u,v)$ so that the parent $z_1$ of $x$ is below $u$, then move the tail of $(z_1,v)$ so that the reticulation $u$ becomes a parent of $x$.}
\label{fig:diamReticMoveEasy}
\end{figure}
 
%	 As the parent of $x$ in $N_1$ is a reticulation, and as $N_1[Y] = N[Y]$,  case~\ref{case:reticulationFound} applies to $N_1$ and $N'$, and so there exists a sequence of $4|N \setminus Y|-4$ tail moves turning $N_1$ into $N'$. As $N_1$ is derived from $N$ by a single tail move, there exists a sequence of $4|N \setminus Y|-3$ tail moves turning $N$ into $N'$.

	 \item{\bf The removal of $\bm{(z,x)}$ followed by suppressing $\bm z$ creates parallel edges:} 
% 	 In this case, let $d$ be the child of $z$ that is not $x$, and let $c$ be the parent of $z$. Then $N$ contains the edges $(c,d),(c,z),(z,d)$. 
	 Then there exists nodes $c \neq x, d \neq x$ such that $c,d,z$ form a triangle with long edge $(c,d)$.
	 As $c$ has outdegree $2$ it is not the root of $N$, so let $b$ denote the parent of $c$.
	 
	    \begin{enumerate}
	    \item\label{case:reticulationMoveTail}  {\bf $\bm b$ is not the root of $\bm N$:}   In this case, let $a$ be a parent of $b$ in $N$. Observe that the edge $(a,b)$ is not below $d$ and that $b \neq d$. Furthermore, $(c,d)$ is movable since $c$ is not a reticulation or the root, and there is no edge 
       % \todo{We need to prove that we cannot have the edge $(b,z)$ and not the edge $(b,d)$, right? I modified the text accordantly} 
        $(b,z)$  (the existence of such an edge would imply that $z$ has indegree $2$, as the edge $(c,z)$ exists. But this is not possible because $z$ is a tree node with child edges $(z,x)$ and $(z,c)$). %$(b,d)$ (such an edge would mean $d$ has indegree $3$, as the edges $(c,d)$ and $(z,d)$ exist).
	    It follows that the tail of $(c,d)$ can be moved to $(a,b)$  (again using Observation~\ref{rem:MovableTo}).	    
	    Let $N_1$ be the network derived from $N$ by applying this tail move (see Figure~\ref{fig:diamReticMoveTail}). 
	    
	    Observe that in $N_1$ we now have the edge $(b,z)$ (as $c$ was suppressed), and still have the edge $(z,d)$ but not the edge $(b,d)$ (such an edge would mean $d$ has indegree $3$, as the edges $(c,d)$ and $(z,d)$ exist). Thus, deleting $(z,x)$ and suppressing $z$ will not create parallel edges, and so $(z,x)$ is movable in $N_1$.	   
	    Thus Case~\ref{case:reticulationMoveEasy} applies to $N_1$ and $N'$, and so there exists a sequence of $3|N \setminus Y|-1$ tail moves turning $N_1$ into $N'$.
	    As $N_1$ is derived from $N$ by a single tail move, there exists a sequence of $3|N \setminus Y|$ tail moves turning $N$ into $N'$.
	    
               \begin{figure}[h!]
\begin{center}
\includegraphics[scale=0.5]{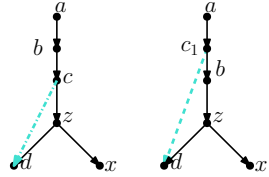}
\end{center}
\caption{
%Illustration of Case~\ref{case:reticulationMoveTail}.
\markj{Proof of Lemma~\ref{lemma:GreenLineLemma}, Case~\ref{case:reticulationMoveTail}:}
If $(z,x)$ is not movable because of the triangle with long edge $(c,d)$, and the parent of $c$ is not the root of $N$, then we move the tail of $(c,d)$ 'further up' in order to make $(z,x)$ movable. }
\label{fig:diamReticMoveTail}
\end{figure}
% Note: In this figure, b is displayed without one of its incident edges, because we don't know or care whether b is a reticulation or a tree node. Is this confusing? Would it be better to add an undirected edge? Add a directed edge and put a disclaimer in the caption saying that the edge could go in the other direction?	   
       \item\label{case:reticulationMoveHead} {\bf $\bm b$ is the root of $\bm N$:} In this case, we observe that 
% 	    $d$ is a reticulation in $N$, and that 
	    if $d \in Y$ then every reticulation in $N$ is in $Y$. This contradicts the fact that $N\setminus Y$ and $N'\setminus Y'$ contain the same number of reticulations. Therefore we may assume that $d \notin Y$.
	    Then we may proceed as follows.
	    Let $N_3$ be the network derived from $Y$ by moving the \emph{head} of $(c,d)$ to $(z,x)$.
        % (see Figure~\ref{fig:diamReticMoveHead}).
       As this is a head move of type (a), it can be replaced with a sequence of four tail moves (see Figure \ref{fig:HeadMoveA'}). However, we show that in this particular case, it is possible to replace it with a sequence of only three tail moves.
       
       Let $e$ be the child of $d$ in $N$. We note that if $e$ is a reticulation, then one of the parents of $e$ is a descendant of $x$ (otherwise $z,b$ or $c$ would have to be a parent of $e$, which is not the case). Thus if $e$ is a reticulation then it is a descendant of $x$, and by a similar argument if $x$ is a reticulation then it is a descendant of $e$.
Thus, we may assume that one of $e,x$ is not a reticulation, and that furthermore at least one of $e,x$ is a tree node (if both are leaves then $N$ and $N'$ are the networks depicted in Figure~\ref{fig:HeadMoveACounter}, and if one is a reticulation then the other must be an ancestor of it, and is therefore a tree node).

Our approach in this case will be to ``swap'' the positions of $e$ and $x$ via a series of tail moves. 
We will assume that $x$ is a tree node, with children $s$ and $t$. 
(The case that $e$ is a tree node can be handled in a similar manner.)
Then we apply the sequence of tail moves depicted in Figure~\ref{fig:diamReticMoveHead}.

%	    We observe that this move cannot create a cycle (such a cycle is only possible if there is a path in $N$  from $x$ to $c$, a contradiction as $x \in Y, c \notin Y$), and cannot create parallel edges (such parallel edges would have to be out-edges of $z$, and this does not happen as one of the out-neighbors of $z$ is the newly created node subdividing $(z,x)$).
%	    Thus, this is a valid head move. 
%	    Furthermore, this is a head move of type (a), and therefore can be replaced with a sequence of $4$ tail moves.
	    
	    Observe that in $N_3$,  $x$ has a parent not in $Y$ which is a reticulation, and that $N_3[Y] = N[Y]$. Then Case~\ref{case:reticulationFound} applies to $N_3$ and $N'$, and so there exists a sequence of $3|N \setminus Y|-3$ tail moves turning $N_3$ into $N'$.
	    As $N_3$ can be derived from $N$ by a sequence of $3$ tail moves, it follows that there exists a sequence of $3|N \setminus Y|$ tail moves turning $N$ into $N'$.
        
                               \begin{figure}[h!]
\begin{center}
\includegraphics[scale=0.5]{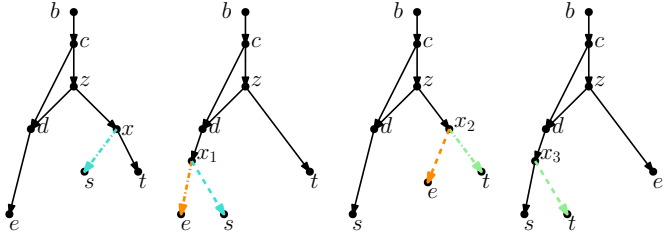}
\end{center}
\caption{
%Illustration of Case~\ref{case:reticulationMoveHead}. 
\markj{Proof of Lemma~\ref{lemma:GreenLineLemma}, Case~\ref{case:reticulationMoveHead}:}
Note that the set of tail moves depicted is equivalent to moving the head of $(c,d)$ to $(z,x)$. After this sequence of moves, $x$ now has a reticulation parent.}
\label{fig:diamReticMoveHead}
\end{figure}
        
 %                      \begin{figure}[h!]
%\begin{center}
%\includegraphics[scale=0.5]{diamReticMoveHead}
%\end{center}
%\caption{Illustration of case~\ref{case:reticulationMoveHead}. If $(z,x)$ is not movable because of the triangle with long edge $(c,d)$, and the parent of $c$ is the root of $N$, then we move the head of $(c,d)$ to $(z,x)$, so that $x$ now has a reticulation parent.}
%\label{fig:diamReticMoveHead}
%\end{figure}

	    \end{enumerate}
	 
	 \end{enumerate}
 
    \end{enumerate}
    
 \item\label{case:reticulationReverse} {\bf There exists a lowest node $\bm u$ of $\bm{N \setminus Y}$ such that $\bm u$ is a reticulation:} By symmetric arguments to Case~\ref{case:reticulation}, we have that there is a sequence of at most $3|N\setminus Y|$ tail moves turning $N'$ into $N$. As all tail moves are reversible, there is also a sequence of at most $3|N\setminus Y|$ tail moves turning $N$ into $N'$.
 
 \item {\bf No lowest node of $\bm{N \setminus Y}$ nor any lowest node of $\bm{N'\setminus Y'}$ is a reticulation:} As $|N \setminus Y| > 1$, we have that in fact every lowest node of $N \setminus Y$ and every lowest node of $N'\setminus Y'$ is a tree node. Then we proceed as follows. Let $u'$ be an arbitrary lowest node of $N' \setminus Y'$, with $x'$ and $y'$ its children.

  Then $x',y'$ are in $Y'$, and therefore there exist nodes $x,y \in Y$ such that $x$ ($y$) is mapped to $x'$ ($y'$) by the isomorphism between $N[Y]$ and $N'[Y']$.
 Furthermore, $x$ has the same number of parents in $N$ as $x'$ does in $N'$, and the same number of parents in $Y$ as $x'$ has in $Y'$ %\textcolor{green}{
 (again because the networks are binary and  $x$ has the same number of children in $N$ as $x'$ does in $N'$ by the isomorphism between $N[Y]$ and $N'[Y']$). %}
 Thus $x$ has at least one parent not in $Y$. Similarly, $y$ has at least one parent not in $Y$.
 
    \begin{enumerate}
     \item\label{case:splitFound} {\bf $\bm x$ and $\bm y$ have a common parent $\bm u$ not in $\bm Y$:} In this case, let $Y_1 = Y \cup \{u\}$ and $Y_1' = Y' \cup \{u'\}$, and extend the isomorphism between $N[Y]$ and $N[Y']$ to an isomorphism between $N[Y_1]$ and $N[Y_1']$, by letting $u$ be mapped to $u'$ (see figure~\ref{fig:diamSplitEasy}).
    We now have that $Y_1$ and $Y_1'$ are downward-closed sets of nodes such that $N[Y_1]$ is isomorphic to $N[Y_1']$, and $L(N)\subseteq Y_1, L(N') \subseteq Y_1'$.
    Furthermore $|N \setminus Y_1| < |N \setminus Y|$.
    Thus by the inductive hypothesis, there is a sequence of $3|N \setminus Y_1| =  3|N \setminus Y|-3$ tail moves turning $N$ into $N'$.
    
\begin{figure}[h!]
\begin{center}
\includegraphics[scale=0.5]{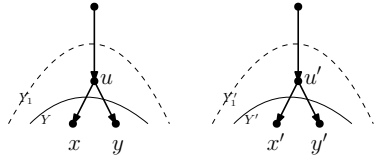}
\end{center}
\caption{
%Illustration of Case~\ref{case:splitFound}.
\markj{Proof of Lemma~\ref{lemma:GreenLineLemma}, Case~\ref{case:splitFound}:}
If $u'$ is a lowest reticulation in $N' \setminus Y'$ with children $x',y'$, and the nodes $x,y \in Y$ corresponding to $x',y'$ share a reticulation parent $u$ in $N \setminus Y$, then we may add $u$ to $Y$ and $u'$ to $Y'$.}
\label{fig:diamSplitEasy}
\end{figure}

    \item {\bf $\bm x$ and $\bm y$ do not have a common parent not in $\bm Y$:} In this case, let $z_x$ be a parent of $x$ not in $Y$, and let $z_y$ be a parent of $y$ not in $Y$. Recall that $z_x$ and $z_y$ are both tree nodes. It follows that either one of $(z_x,x), (z_y,y)$ is movable, or deleting $(z_x,x)$ and suppressing $x$  (deleting $(z_y,y)$ and suppressing $y$) would create parallel edges.
    
	\begin{enumerate}
	  \item\label{case:splitMoveEasy}  {\bf $\bm{(z_x,x)}$ is movable:}  in this case, observe that the edge $(z_y,y)$ is not below $x$ (as $x \in Y$ and $z_y \notin Y$), and that $x \neq y$.  Then by Observation~\ref{rem:MovableTo}, the tail of $(z_x,x)$ can be moved to $(z_y,y)$.
	  
	 Let $N_1$ be the network derived from $N$ by applying this tail move (see Figure~\ref{fig:diamSplitMoveEasy}).
	 Then as $x$ and $y$ have a common parent in $N_1$ not in $Y$, and as $N_1[Y] = N[Y]$, we may apply the arguments of Case~\ref{case:splitFound}  to show that there exists a sequence of $3|N \setminus Y|-3$ tail moves turning $N_1$ into $N'$. As $N_1$ is derived from $N$ by a single tail move, there exists a sequence of $3|N \setminus Y|-2$ tail moves turning $N$ into $N'$.
	 
     \begin{figure}[h!]
\begin{center}
\includegraphics[scale=0.5]{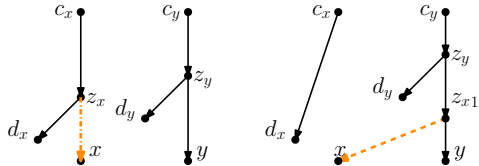}
\end{center}
\caption{
% Illustration of Case~\ref{case:splitMoveEasy}.
\markj{Proof of Lemma~\ref{lemma:GreenLineLemma}, Case~\ref{case:splitMoveEasy}:}
If $(z_x,x)$ is movable, we may move the tail of $(z_x,x)$ to $(z_y,y)$ so that $x$ and $y$ share a parent. }
\label{fig:diamSplitMoveEasy}
\end{figure}
     
	 \item\label{case:splitMoveReverse}  {\bf $\bm{(z_y,y)}$ is movable:} By symmetric arguments to Case~\ref{case:splitMoveEasy}, we have that there is a sequence of at most  $3|N \setminus Y|-2$  tail moves turning $N$ into $N'$.
	 
	 \item\label{case:splitMoveTriangles} {\bf Neither $\bm{(z_x,x)}$ nor $\bm{(z_y,y)}$ is movable:} 
% 	 Hoo boy, we're getting pretty deep into it now. How you holding up out there? Don't lose hope; we've just got to get through this last big push, and then we're on the home straight.
	 In this case, there must exist nodes $d_x,c_x, d_y, c_y$ such that 
% 	 edges $(z_x,d_x),(c_x,z_x),(c_x,d_x)$ and $(z_y,d_y),(c_y,z_y),(c_y,d_y)$ exist.
	 $c_x,d_x,z_x$ form a triangle with long edge $(c_x,d_x)$, and  $c_y,d_y,z_y$ form a triangle with long edge $(c_y,d_y)$.
	 Moreover, as $z_x,z_y$ are different nodes with one parent each, $c_x \neq c_y$. It follows that one of $c_x,c_y$ is not the child of the root of $N$. Suppose without loss of generality that $c_x$ is not the child of the root. Then there exist nodes $a_x,b_x$ and edges $(a_x,b_x), (b_x,c_x)$.
	 
	 By similar arguments to those used in Case~\ref{case:reticulationMoveTail}, the tail of $(c_x,d_x)$ can be moved to $(a_x,b_x)$, and in the resulting network $N_1$, $(z_x,x)$ is movable (see Figure~\ref{fig:diamSplitMoveTail}).
	 Thus Case~\ref{case:splitMoveEasy} applies to $N_1$ and $N'$, and so there exists a sequence of $3|N\setminus Y|-2$ tail moves turning $N_1$ into $N'$. As $N_1$ is derived from $N$ by a single tail move, there exists a sequence of $3|N\setminus Y|-1$ tail moves turning $N$ into $N'$.\qed
     
          \begin{figure}[h!]
\begin{center}
\includegraphics[scale=0.5]{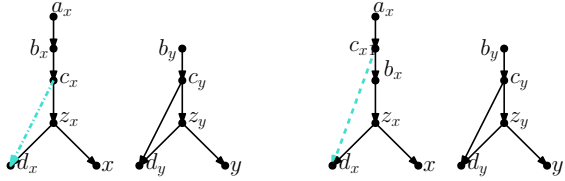}
\end{center}
\caption{
%Illustration of Case~\ref{case:splitMoveTriangles}.
\markj{Proof of Lemma~\ref{lemma:GreenLineLemma}, Case~\ref{case:splitMoveTriangles}:}
If neither $(z_x,x)$ or $(z_y,y)$ is movable, we can make at least one of them movable by moving the long edge of its triangle ``further up''.}
\label{fig:diamSplitMoveTail}
\end{figure}
	\end{enumerate}
    \end{enumerate}
 \end{enumerate}
\end{proof}

 By setting $Y = L(N)$ and $Y' = L(N')$, we have the following:%\todo{Leo: We need to say something about rSPR here, because we use it below.}
 \begin{theorem}\label{the:GreenLine}
 Let $N$ and $N'$ be networks
 %\textcolor{green}{
 in the $k$-th tier on $X$
 %}
 %on $X$ with the same number of reticulations,
such that $N$ and $N'$ are not the networks depicted in Figure~\ref{fig:HeadMoveACounter}. 
 Then there is a sequence of at most $3(|N|-|X|)=3(|X|+2k)$ tail moves turning $N$ into $N'$.
\end{theorem}

%By similar arguments, we can get a bound on the number of rSPR moves needed to turn $N$ into $N'$.

 As rSPR moves consist of head moves and tail moves, Theorem~\ref{the:GreenLine} also gives us an upper bound on the number of rSPR moves needed to turn $N$ into $N'$.
 By modifying these arguments slightly, we can improve this bound in the case of rSPR moves.

\begin{theorem}\label{the:GreenLineSPR}
 Let $N$ and $N'$ be networks in the $k$-th tier on $X$.
 Then there is a sequence of at most $2|X|+3k-1$ rSPR moves turning $N$ into $N'$.
\end{theorem}
\begin{proof}
%In this proof, we will use terminology of the proof of Lemma~\ref{lemma:GreenLineLemma}, all cases referred to are in this proof.\\
%Note that Cases~\ref{case:reticulation} and~\ref{case:reticulationReverse} in the proof of Lemma~\ref{lemma:GreenLineLemma} can both be resolved with one head move. \remie{There is a reticulation node in $N\setminus Y$. Using an argument similar to Remark~\ref{rem:MovableTree}, one of its incoming edges is head-movable. And, because $Y$ is downward-closed, head-moving this edge to $(u,x)$ will not create any cycles or parallel edges. Except if the only movable edge is $(u,v)$, but then the child of $v$ is also a reticulation, and both incoming edges are head-movable, and they can be moved to $(u,x)$.}

Recall that the proof of Lemma~\ref{lemma:GreenLineLemma} works by gradually expanding two downward-closed subsets $Y \subseteq V(N), Y' \subseteq V(N')$ for which $N[Y]$ is isomorphic to $N'[Y']$, using at most $3$ tail moves each time the size of $Y$ and $Y'$ is increased.
 We show that in Cases~\ref{case:reticulation}  and~\ref{case:reticulationReverse} of the proof of Lemma~\ref{lemma:GreenLineLemma}, we may instead use one head move. 
 Indeed, in Case~\ref{case:reticulation} there is a lowest node $u'$ of $N'\setminus Y'$ that is a reticulation with child $x' \in  Y$, 
 and the node $x \in Y$ corresponding to $x'$ has parent $z$. If $z$ is a reticulation (Case~\ref{case:reticulationFound}) then as before there is no need for any move, we simply add $z$ to $Y$ and $u'$ to $Y'$.
 If $z$ is not a reticulation, we proceed as follows. 
 There exists some reticulation node in $v \in N\setminus Y$ (again, such a node must exist, as $u'$ exists and $N\setminus Y$ and $N'\setminus Y'$  have the same number of reticulations). 
 Moving one of its parent edges $(u,v)$ to $(z,x)$ will not create a cycle, as $Y$ is downward-closed.
 It cannot create any parallel edges, unless either $u = z$, or the other parent edge $(w,v)$ is part of a triangle on $w,v$ and the child of $v$.
If $u = z$, then we can move $(w,v)$ to $(z,x)$ and this will not create parallel edges
 %
 %Thus, one of the parent edges of $v$ can be head-moved to $(z,x)$, 
 %unless $(z,v)$ is an edge and 
 unless $z,v,x$ form a triangle. But in this case $x$ is a child of $v$ and thus $x$ already has a parent in $N \setminus Y$ that is a reticulation.
 This implies that, after at most one head move, $x$ has a reticulation parent, and we may proceed as in Case~\ref{case:reticulationFound}.
 (Case~\ref{case:reticulationReverse} is handled symmetrically.)

The other cases use at most one tail move (and thus at most one rSPR move), apart from Case~\ref{case:splitMoveTriangles} that may require 2 rSPR moves. This case can come up as many times as there are tree nodes in the network. 
Hence the number of moves needed to add a node to 
% the downward-closed set 
$Y$ is at most one for each reticulation node, and at most two for each tree node. This means at most $|V|-|X|+t$ moves are needed, where $t$ denotes the number of tree nodes in $N$ (and thus in $N'$).\\
Recall from Observation~\ref{obs:numbers} that $|V|=2(|X|+k)$ for any binary tier-$k$ network.
As $k$ nodes are reticulations, $|X|$ are leaves and $1$ is the root, we have $t = |X|+k-1$.
%and that therefore $t=|X|+k$. 
This shows that we need at most $2(|X|+k)-|X|+|X|+k-1=2|X|+3k-1$ rSPR moves to turn $N$ into $N'$.\qed
\end{proof}

In practice, we expect the distance between most pairs of networks to be less than $|X|+3k-1$ because only one case needs two rSPR moves, and this case might not come up very often. %However, it is impossible to exclude this case, this can be shown by example.}
%\todo{RJ: we may also remove this, it doesn't add much, it only indicates that I think that $\Delta^{\text{rSPR}}_k=n+2k+o(n)$ or something like this.}

%\subsection{Diameters compared}
%%\todoCeline{Tail with or without capital letter? It was without in the previous sections, and with now...}
%Tail moves are related to other types of moves such as rNNI and rSPR, which has implications for their induced metrics and diameters. There are obvious bounds when a class of moves contains another class of moves (e.g. $\text{rNNI}\subseteq\text{rSPR}$), by which we mean that each move of the second class is equivalent to a move of the first class. 
%Denote by $\Delta_k^{\mathcal{M}}(n)$ the \emph{diameter} of the $k$-th tier  with $n$ leaves with respect to $\mathcal{M}$; that is, the maximum value of $d_M(N,N')$ over any pair of networks $N,N'$ belonging to the $k$-th tier on a fixed leaf set $X$ with $|X|=n$. Then
%% Denoting with $\Delta_k^{\mathcal{M}}(n)$  the \emph{diameter} of tier $k$ with $n$ leaves and metric induced by move $\mathcal{M}$, 
%we have that $\mathcal{M}\subseteq \mathcal{M}'$ implies $\Delta_k^{\mathcal{M}}(n)\geq \Delta_k^{\mathcal{M}'}(n)$.%Because any $M$-sequence is a $M'$-sequence.

\celine{The next \markj{observation} will be useful \markj{for obtaining} lower bounds for the diameter of tail and rSPR moves.}
\begin{observation}
\label{remark:impliedBounds}
Let $N$ and $N'$ be networks in the $k$-th tier on $X$. The observation that $\text{Tail}_1\subseteq\text{rNNI}\subseteq\text{rSPR}$ %\citep{Gambette-vI-Rearrangement} \todo{Do we actually observe that $\text{Tail}_1\subseteq\text{rNNI}$ anywhere? \celine{Celine: I suggested to discuss this in section 2.2}} 
(where $\text{Tail}_1$ denotes the class of distance-1 tail moves) implies that
\[d_{\text{Tail}_1}(N,N')\geq d_{\text{rNNI}}(N,N')\geq d_{\text{rSPR}}(N,N'),\]
and 
\[\Delta_k^{\text{Tail}_1}(n)\geq\Delta_k^{\text{rNNI}}(n)\geq\Delta_k^{\text{rSPR}}(n).\]
Similarly, $\text{Tail}_1\subseteq\text{Tail}\subseteq\text{rSPR}$ implies that
\[d_{\text{Tail}_1}(N,N')\geq d_{\text{Tail}}(N,N')\geq d_{\text{rSPR}}(N,N'),\]
and
\[\Delta_k^{\text{Tail}_1}(n)\geq\Delta_k^{\text{Tail}}(n)\geq\Delta_k^{\text{rSPR}}(n).\]\qed
\end{observation}

Diameters of move-induced metrics on tree space are well studied. A few relevant bounds are
$\Delta_0^{\text{rSPR}}=n-\Theta(\sqrt{n})$ \citep{ding2011agreement,atkins2015extremal} and $\Delta_0^{\text{rNNI}}=\Theta(n\log(n))$ \citep{li1996nearest}. We extend these results to higher tiers of network space.

Lemma~\ref{lem:distanceRelationMycorrhizal} %\todo{move this lemma back to here? Leo: yes! If we use it here we should move it here. \celine{Celine: I agree. Maybe move section 4.4 before section 4.2 if you want to keep the sections balanced ?}} 
gives lower bounds of order $n-\Omega(\sqrt{n})$ on diameters for \leo{rSPR} and tail moves, by reducing to trees and using the corresponding diameter bound. Theorem~\ref{the:GreenLine} and Theorem~\ref{the:GreenLineSPR} give upper bounds of order $n$. More precisely, $\Delta_k^{\text{Tail}}(n)\leq 3(n+2k)$ and $\Delta_k^{\text{rSPR}}(n)\leq 3n+2k$ from Observation~\ref{remark:impliedBounds}. The following theorem summarizes this discussion:  
\begin{theorem}
The diameter of tiers of network space for metrics induced by rSPR and tail moves satisfy
\begin{equation*}
\begin{aligned}
n-\Theta(\sqrt{n})&\leq&\Delta_k^{\text{Tail}}(n)&\leq& 3(n+2k),\\
n-\Theta(\sqrt{n})&\leq&\Delta_k^{\text{rSPR}}(n)&\leq& 2n+3k.
\end{aligned}
\end{equation*}
\end{theorem}

Additionally, going back through the proof of Theorem~\ref{the:HeadToTail}, we see that any head move can be replaced by at most four tail moves if the network has more than one leaf, so $d_{\text{Tail}}(N,N')\leq 4 d_{rSPR}(N,N')$ for any pair of networks $N,N'$ with more than one leaf. %If the networks have one leaf, such a result is unavailable to us, as we did not find a way to rewrite some of the head moves into short sequences of tail moves. 
Note that this does not give us new information regarding the diameters, but it does give bounds on the distances when we are given two networks.
%\todo{Say something about sharpness of bounds?}
%
\celine{\subsection{The diameter of Tail\textsubscript{1} and rNNI moves}}

%\subsection{Tail{\textsubscript{1}}{ distance}}
Comparing rNNI and tail moves directly is complicated by the fact that rNNI moves are more local. To make the comparison easier, we consider local tail moves: tail moves over small distance. The following lemma indicates how restricting to distance-1 tail moves influences our results.

\begin{lemma}\label{D1TailSub}
Let $e=(u,v)$ to $f=(s,t)$ be a valid tail move in a tier $k$ network $N$ on $X$, then there is a sequence of at most $|X|+3k-1$ distance-1 tail moves resulting in the same network.
\end{lemma}
\begin{proof}
Note that there exist directed paths $\LCA(u,s)\rightarrow~s$ and $\LCA(u,s)\rightarrow~u$ (which are not necessarily unique) for any choice of $\LCA(u,s)$. We prove that a tail move of $e$ to any edge on either path is valid, and this gives a sequence of distance-1 tail moves: Indeed, for any edges $f$ and $g$ that share a node, if there is a valid tail move $e$ to $f$ resulting in network $N_f$, and a valid tail move $e$ to $g$ resulting in network $N_g$, then there is a distance-1 tail move $f$ to $g$ that transforms $N_f$ into $N_g$, and furthermore as $N_g$ has no cycles or parallel edges, this is a valid tail move.

Let $g=(x,y)$ be an edge of one of these two paths. We first use a proof by contradiction to show that the move to $g$ cannot create cycles, then we prove that we do not create parallel edges.

Suppose moving $e$ to $g$ creates a cycle, then this cycle must involve the new edge $e'$ from $g$ to $v$. This means there is a path from $v$ to $x$. However, $x$ is above $u$ or above $s$, which means that either the starting graph is not a phylogenetic network, or the move $e$ to $f$ is not valid. From this contradiction, we conclude that the move of $e$ to $g$ does not create cycles.

Note that $e$ is movable, because the move of $e$ to $f$ is valid. Hence the only way to create parallel edges, is by moving $e=(u,v)$ to an edge $g=(x,y)$ with $y=v$. It is clear that $g$ is not in the path $\LCA(u,s)\rightarrow u$, as this would imply the existence of a cycle in the original network. Hence $g$ must be on the other path. If $f=g$, then the original move of $e$ to $f$ would create parallel edges, and if $g$ is above $f$, the original move moves $e$ to below $e$ creating a cycle. We conclude that there cannot be an edge $g=(x,y)$ on either path such that $y=v$, hence we do not create parallel edges.
%Can be done shorter and more easily if I first do the lemma proving that moving up is allowed if an edge is movable.
%Only be careful of isomorphic networks, even though they are allowed, they will never help us.

Noting that a path between two nodes uses at most $|E|-|X|$ edges, %\todo{Explain this more? The path involves at most either 2 leaves and no root, or one leaf and the root; so at most 2 pendant edges, of which we count only one in the path length, and there are $|X\cup\{\rho\}|$ pendant edges.}
we see that we need at most $|E|-|X|=|X|+3k-1$ distance-1 tail moves to simulate a long distance tail move (the last equivalence holds by Observation \ref{obs:numbers}).
\qed\end{proof}

Note that a distance-$d$ tail move cannot necessarily be simulated with a sequence of $d$ distance-1 tail moves. The path of length $d$ defining the distance of the tail move might not be a path over which we can move the tail: it might go down to a reticulation, and then up again.

Lemma~\ref{D1TailSub} directly gives us upper bounds on Tail$_1$ and rNNI diameters in terms of the tail diameter: each tail move is replaced by distance-1 tail moves, giving an upper bound of $(3n+6k)(n+3k-1)$ for tier $k$ networks on $n$ leaves. As we are mostly interested in the effect of the number of leaves, we denote these bounds $\Delta_k^{\text{Tail}_1}(n)=O(n^2)$ and $\Delta_k^{\text{rNNI}}(n)=O(n^2)$ (because $\text{Tail}_1\subseteq\text{rNNI}$).

%Francis et al. \citep{francis2017bounds} 
\markj{Francis et al. (\citeyear{francis2017bounds})}
proved diameter bounds for $d_{\text{NNI}}$ on \emph{unrooted} networks. These moves do not have to account for the orientation of edges. Therefore they generally define larger classes of moves.
%, e.g. ``$\text{rNNI}\subseteq \text{NNI}$''.
More precisely, given a rooted network $N$ with unrooted underlying graph $U(N)$, the set of unrooted networks that can be reached with one NNI move from $U(N)$ contains the set of unrooted networks we get by applying one rNNI move to $N$ and then taking the underlying graph. Francis et al. give the following lower bound on NNI diameters using Echidna graphs~\citep[Theorem 4.3]{francis2017bounds}:
\[\left[(v_k^n-3)\log_6\left(\frac{v_k^n}{2}-2\right)-(2k-1)\log_6(k-1)-(v_k^n-2k)\log_6 e -2v_k^n\right],\]%lower bound rooted NNI diameter.
where $v_k^n$ is the number of nodes in an unrooted network with $n+1$ leaves and $k$ reticulation nodes. This lower bound is $\Omega(n \log(n))$ for fixed $k$.
As Echidna graphs are rootable (i.e., for each Echidna graph there exists some rooted network with this Echidna graph as underlying unrooted network), the argument of Francis et al. easily extends to rooted networks and we get the same lower bound for $\Delta_k^{\text{rNNI}}(n)$.
The preceding discussion proves the following theorem:

\begin{theorem}
The diameter of tiers of network space for metrics induced by rNNI and distance-1 tail moves satisfy 
\[\Delta_k^{\text{rNNI}}(n)=\Omega(n \log (n))\quad\text{ and }\quad\Delta_k^{\text{rNNI}}(n)=O(n^2),\]
\[\Delta_k^{\text{Tail}_1}(n)=\Omega(n \log (n))\quad\text{ and }\quad \Delta_k^{\text{Tail}_1}(n)=O(n^2).\]
\end{theorem}
%\todo{Do we want this to be more precise as for tail and rSPR moves? It feels less useful, because this will probably be improved to $\Theta(n\log(n))$ anyway.}
%\todo{Can we alternatively prove $\Delta_k^{\text{NNI}}\leq\Delta_k^{\text{rNNI}}$? This might not be that easy: some unrooted networks do not correspond to rooted networks.}
%\todoCeline{Answer to the 2 previous questions: I would keep it easy. The paper is getting longer...}
\subsection{The diameter of SPR moves}\label{sec:unrooted}
%\todoCeline{The new stuff is here}

\celine{In this subsection, we will give an upper bound on the diameter of SPR moves on unrooted networks, using results for  rSPR moves on rooted networks.}

The underlying unrooted network of a rooted network $N$ is denoted $U(N)$.
%\begin{definition}[Rootable network]%take care: the root should be chosen to be one of the leaves!
An unrooted network $U$ is called \emph{rootable} if there exists a rooted network $N$ such that $U(N)=U$.  
%\end{definition}
%
As 
% already 
mentioned at the end of Section \ref{sec:rearrangement}, not all unrooted networks are rootable \markj{(see Figure \ref{fig:HuisjeBoompjeBeestje})}.

Any edge of an unrooted network whose removal disconnects the network is called a \emph{cut-edge}. If only one of the components contains leaves, the edge is called a \emph{redundant cut-edge}. \leo{A \emph{blob} of an unrooted network is a nontrivial biconnected component, i.e. a maximal subgraph with at least two vertices and no cut-edges.}
\celine{The next lemma characterizes rootable networks via redundant cut-edges:}

%\todoRemie{Leo (or anyone else), can you tell me if this is a known result? That would save some paper in the paper. Leo: I haven't seen this anywhere else.}
\begin{lemma}\label{lem:redundantcutedge}
\leo{An unrooted network
% leo: actually, unrooted networks have at least two leaves by definition
is rootable \markj{if and only if} it has no redundant cut-edges.}
\end{lemma}
\begin{proof}
\leo{First let~$U$ be some unrooted network with no redundant cut-edges. To show that~$U$ is rootable}, we pick any leaf $r$ of $U$ and show how to construct a rooted network with $U$ as underlying graph and $r$ as root. 
First, \leo{orient all cut-edges ``away'' from~$r$.}
%note that any unrooted network can be viewed as a blobbed-tree \citep{gusfield2005fundamental}, 
%\leo{observe that contracting each blob into a single node gives a tree}. %Correctly orienting the edges of $U$ that are in its blobbed-tree is simple: for each edge $\{s,t\}$,  we orient the edge from $s$ to $t$ if there is a path from $s$ to $r$ avoiding $t$ and, otherwise, we orient it from $t$ to $s$.
%Clearly, all edges of the \leo{tree} can be oriented ``away'' from $r$.
Then, it only remains to find a valid orientation of every \leo{blob}.
To this end, let $B$ be \leo{a blob}.
After orienting \leo{all cut-edges}, $B$ has only one incoming edge $(s,t)$ and, as $(s,t)$ is not redundant, $B$ has at least one outgoing edge $(x,y)$.
Since $B$ is biconnected, there is a \leo{bipolar (i.e. acyclic)} orientation of $B$ with $t$ as source and $x$ as sink \citep{lempel1967algorithm}.
%This gives us an acyclic orientation of $B$ directed away from $r$  and in the direction of at least a leaf (the edge $(x,y)$ is not a redundant cut-edge so there is at least one leaf under it), and thus, 
Doing the same for all biconnected components, we get an acyclic orientation of $U$ rooted at $r$.

Conversely, suppose that $U$ has a redundant cut-edge\leo{~$e$. Deleting~$e$ creates a component~$C$ without leaves. If we direct~$e$ towards~$C$, then~$C$ has one source and no possible sinks (no leaves). If we direct~$e$ away from~$C$ then~$C$ has one sink but no possible sources (since the root is also a leaf).} This implies there is no valid orientation of the edges \leo{in~$C$} and therefore in $U$.\qed
\end{proof}

% \begin{observation}
% \todo{RJ: this might be useful for an even better bound. Celine: if we use it, we should rewrite it in a clearer way
% }
% From the proof of the previous lemma it also follows that, for each blob, we can choose any outgoing arc as the sink, meaning that the top vertex of this edge is a reticulation.
% \end{observation}

A \emph{\leo{redundant} terminal component} of a network $U$ is a \leo{nontrivial} biconnected component \leo{that is incident to exactly one cut-edge (which must be a redundant cut-edge).}
% $B$ which corresponds to a leaf in the blobbed tree of $U$
% such that $B$ does not contain any leaf of $U$.
% \celine{A \emph{terminal component} of a network $U$ is a biconnected component $B$ for which there exists only one edge in $U$, say  $\{s,t\}$,  with $s\in B$  and  $t\not\in B$.
% A \emph{non-labeled terminal component} is a terminal component which does not contain leaves.}
%
\leo{The next lemma, which follows directly from Lemma~\ref{lem:redundantcutedge}, characterizes rootable networks via redundant terminal components.}
\begin{lemma}\label{lem:rootableIfNoNLTC}
\leo{A network is rootable if and only if it has no redundant terminal components.}
\end{lemma}
% LEO: I don't think we need a proof for this.
%
% \begin{proof}
% \celine{We shall prove that a network without \leo{redundant} components has no redundant cut-edges. The result will then follow by the previous lemma.}
% Let $U$ be a network with a redundant cut-edge, then one of the components induced by its removal has no leaves. This component contains biconnected components, at least one of which must be a leaf in the blobbed tree as the network is finite. % \celine{and binary}.\todoRemie{Why does it need to be binary? The redundant cut-edge represents a cut-edge in the blobbed tree, even if $U$ is non-binary. On the non-leaf side of this cut edge, there must be a leaf (as it is a tree). This leaf can only be a blob, as there are no labels of $U$ there.} 
% Hence this component is a non-labeled terminal component. We conclude that a network without non-labeled terminal components cannot contain any redundant cut-edge.\qed
% \end{proof}

We now give a formal definition of a SPR move: 
\begin{definition}
Let $U$ be an unrooted network, and let $\{u,v\}$, $\{x,u\}$, $\{u,y\}$ and $e$ be edges of $U$.   
The \emph{SPR move} of the $u$-end of $\{u,v\}$ to $e$ consists of the subdivision of $e$ with $u'$, the removal of $\{u,v\}$, the suppression of $u$, and the addition of edge $\{u',v\}$. The move is only valid if the resulting graph is an unrooted network, i.e. if it is connected and has no parallel edges.
\end{definition}

\celine{The next two lemmas give upper bounds for the number of \leo{redundant} terminal components in an unrooted network with reticulation number $k=|E|-|V|+1$, and for the number of SPR moves needed to get to a network without redundant terminal components.}
\begin{lemma}\label{lem:toRootable}
Let $U$ be an unrooted network in the $k$-th tier. Then $U$ has at most $k/3$ \leo{redundant} terminal components.
\end{lemma}
\begin{proof}
Let $V'$ be the nodes in a \leo{redundant} terminal component, and $E'$ the edges with at least one \leo{endpoint} in $E'$. Then every node in $V'$ has degree $3$ \leo{w.r.t.~$E'$} and every edge in $E'$ except one has both \leo{endpoints} in $V'$, which implies that $3|V'| = 2|E'|- 1$. It follows that $|V'|$ is odd. Furthermore,~$V'$ must be greater than~$3$ in order for every node to have degree $3$. Thus, $|V'|\geq 5$, and hence $|E'|-|V'| = (3|V'|+1)/2 - |V'| = (|V'|+1)/2 \geq (5+1)/2 = 3$. 
Now observe that every node and edge \leo{of the network} appears in at most one such set $V'$ or $E'$. Furthermore, if all such nodes and edges are deleted, the resulting graph $(V'',E'')$ is still connected and so satisfies $|E''| - |V''| \geq -1$.
It follows that if $U$ has more than $k/3$ \leo{redundant} terminal components, then the reticulation number of $U$ is $|E|-|V|+1 > 3(k/3)+|E''|-|V''| + 1\geq k$, a contradiction.\qed

%A non-labeled terminal component contributes at least $6$ to the number of nodes  of $U$: \remie{$4$ for any cubic graph, and $2$ to connect the cubic graph to the network (for a minimal example of a non-labeled terminal component, see Figure~\ref{fig:HuisjeBoompjeBeestje}). All of these nodes are not part of the tree structure on the leaves.}  \todoCeline{You should prove that you cannot contribute 0 (obvious) or 1}\todoRemie{I added a definition of reticulation number somewhere above, is that enough to explain this? It basically says that each reticulation adds two nodes.}
%For each two extra nodes, we need exactly one reticulation, thus, a network with $c$ non-labeled terminal components is in tier greater or equal to $3c$. This implies that a network in the $k$-th tier has at most $k/3$ non-labeled terminal components.\qed
\end{proof}

\begin{lemma}\label{lem:reduceNLTC}
Let $N$ be an unrooted network in the $k$-th tier with $c$ \leo{redundant} terminal components, then there exists an unrooted  network $N'$ in the $k$-th tier with \leo{at most} $c-1$ \leo{redundant} terminal components such that $d_{\textrm{SPR}}(N,N')=1$.
\end{lemma}
\begin{proof}
Pick any \leo{redundant} terminal component $B$ and let $\{u,v\}$ be 
%the redundant cut-edge for this component, with the non-labeled terminal component connected to $v$. 
the unique edge for which $u \notin B$, $v \in B$.
Let $x$ and $y$ be the other neighbours of $v$. Now SPR move the $v$-end of edge $\{v,x\}$ to a leaf edge of $U$. Suppressing $v$ cannot give parallel edges, because $\{u,v\}$ is a cut-edge (Figure~\ref{fig:TerminalBlobRemoval}). In the resulting network, $B$ is extended to a biconnected component with a pendant leaf, and because no new cut-edges have been created, the network has \leo{at most} $c-1$ \leo{redundant} terminal components and is one SPR move away from the original network. Note that the networks are in the same tier because \leo{an} SPR move does not change the number of 
edges 
%leaves 
nor the number of vertices.\qed

\begin{figure}[h!]
\begin{center}
\includegraphics[scale=0.5]{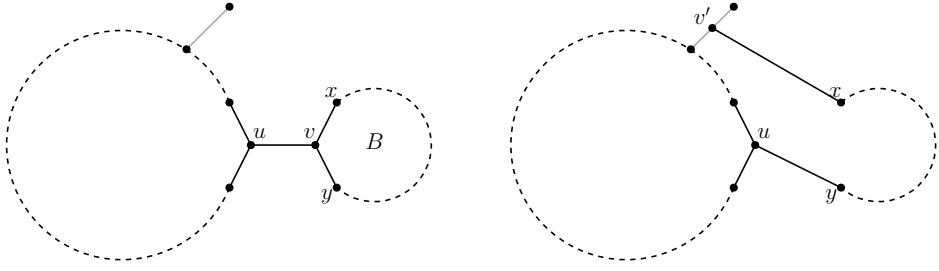}
\end{center}
\caption{The SPR move used to remove \leo{redundant} terminal component $B$. The big circle represents all of $U$ except for $B$. Note that in the network on the right, $B$ is not a \leo{redundant} terminal component anymore, and no extra \leo{redundant} terminal components have been created.}
\label{fig:TerminalBlobRemoval}
\end{figure}
\end{proof}

\markj{Lemmas~\ref{lem:rootableIfNoNLTC},~\ref{lem:toRootable} and~\ref{lem:reduceNLTC} imply the following:}
\begin{corollary}\label{cor:ToRootable}
Any tier $k$ network is at most $k/3$ SPR moves away from a rootable network.
\end{corollary}

\celine{The next two results 
%make a connection between 
\markj{relate}
a set of rSPR moves between two rooted networks and the corresponding set of SPR moves  between their underlying graphs.}
\begin{lemma}
Any rSPR move transforming a rooted network $N$ into another \leo{rooted network} $N'$ \leo{has a corresponding} SPR move transforming $U(N)$ into $U(N')$.
\end{lemma}
\begin{proof}
Suppose that \leo{the rSPR move} is a head move of $(u,v)$ to $e$. \leo{Then} $U(N')$ is the network we get by doing the SPR move of the $v$-end of $\{u,v\}$ to $e$. This SPR move is valid because $N'$, and therefore $U(N')$ is connected with no parallel edges. \leo{Similarly, a tail move of~$(u,v)$ to~$e$ has a corresponding SPR move moving the $u$-end of~$\{u,v\}$ to~$e$.}\qed
\end{proof}

\begin{corollary}\label{cor:rootableSPRsequences}
%Any rSPR sequence of rooted networks corresponds to an SPR sequence of rootable networks. Hence the diameter for rooted networks gives an upper bound for the diameter of rootable networks (unrooted).
Any rSPR sequence transforming a rooted network $N$ into another \leo{rooted network} $N'$ \leo{has a corresponding} SPR sequence transforming  $U(N)$ into $U(N')$. 
\end{corollary}
%\celine{Since  $U(N)$ and $U(N')$ are rootable networks,}
\leo{Moreover, since any rootable network is \celine{the underlying graph} $U(N)$ of some rooted network $N$, 
the diameter for rSPR on rooted networks gives an upper bound for the diameter of SPR moves on unrooted rootable networks.} %LEO: I took this out of the corollary

\leo{We do not apply the following proposition here, but it may become useful when a better bound on $\Delta_k^{\textrm{rSPR}}(n)$ is found in future research.}

\begin{proposition}
%\todoRemie{This proposition is not necessary, but if a better bound on $\Delta_k^{\textrm{rSPR}}(n)$ is found, it might become useful again.}
The diameter of the $k$-tier  of unrooted network space has upper bound
\[\Delta_k^{\textrm{SPR}}(n)\leq \Delta_k^{\textrm{rSPR}}(n)+2M,\]
where $M$ denotes the maximal SPR distance from any unrooted network to a rootable network.
\end{proposition}
\begin{proof}
Let $N$ and $N'$ be tier $k$ unrooted networks. Then we can do $M$ moves from $N$ and $M$ moves from $N'$ to get to rootable networks $N_r$ and $N'_r$. Choose a root for both these networks, then by Corollary~\ref{cor:rootableSPRsequences} there is a sequence of at most $\Delta_k^{\textrm{rSPR}}(n)$ moves to go from $N_r$ to $N'_r$. Because all moves are reversible, there is a sequence of moves:

\[N\xleftrightarrow{M} N_r \xleftrightarrow{\Delta_k^{\textrm{rSPR}}(n)} N'_r \xleftrightarrow{M} N',\]
of length at most $\Delta_k^{\textrm{rSPR}}(n)+2M$ from $N$ to $N'$.\qed

\end{proof}

This proposition together with Corollary~\ref{cor:ToRootable} and Theorem~\ref{the:GreenLineSPR} gives a reasonable bound on the diameter of unrooted networks, namely $2|X|+3k+\frac{2}{3}k$. However, using the (lack of) structure in an unrooted network, we can do better. The next theorem again uses rooted networks, but dynamically re-orientates the network during the induction.

\begin{theorem}\label{thm:proofNNIimproved}
The \leo{SPR} diameter of the $k$-tier  of unrooted network space has upper bound
\[\Delta_k^{\textrm{SPR}}(n) \leq n+\leo{\frac{8}{3}k}.\]
\end{theorem}
\begin{proof}
Let $N$ and $N'$ be tier $k$ networks. \leo{As before}, we use $\frac{2}{3}k$ moves to produce rootable networks $N_r$ and $N'_r$. Now we do induction as in the proof of the rooted rSPR diameter \leo{(Theorem~\ref{the:GreenLineSPR})}. Choose a network orientation on both networks. As before, let $Y$ and $Y'$ be the downward-closed subsets of $N_r$ and $N'_r$ making $N_r|_{Y}$ and $N'_r|_{Y'}$ isomorphic. We prove that we need at most $n+2k$ moves in total to produce a network $N$ from $N_r$ and $N'_r$, by inductively increasing the size of $Y$ and $Y'$.

As in the proof \leo{of Theorem~\ref{the:GreenLineSPR}}, we can add a node to $Y$ and $Y'$ using at most one move, except in the case that all lowest nodes are tree nodes and the nodes corresponding to their child nodes have unmovable incoming edges. % Suppose there exists a lowest node $u$ in $N_r\setminus Y$ such that $u$ is a reticulation. Let $x$ be the child of $u$ and note that $x$ corresponds to some node $x'$ in $Y'$. As in the rooted case we need to do at most one move in $N_r'$ to make sure that $x'$ is also directly below a reticulation, so we can add one node to $Y$ and $Y'$ using at most one move.
% By a symmetric argument, if there exists a lowest node $u'$ in $N_r'\setminus Y'$ such that $u'$ is a reticulation, we can add one node to $Y$ and $Y'$ using at most one move.
% Now suppose all lowest nodes above $Y$ and $Y'$ are tree nodes. Let $u$ be such a node with children $x,y\in Y$. Then $x$ and $y$ correspond to some $x'$ and $y'$ in $Y'$. As proof of rooted case provides: if one of the parent edges of $x'$ or $y'$ is movable, we need at most one move to add $u$ and a corresponding node to $Y$ and $Y'$. If neither of the parent edges are movable, the rooted case needs two moves to add a node to $Y$ and $Y'$. 
We now show how this case can be treated differently if we consider SPR moves on unrooted networks.

Suppose we are in the situation of Case~\ref{case:splitMoveTriangles}. We shortly recall the situation:
\emph{Every lowest node of $N \setminus Y$ and every lowest node of $N'\setminus Y'$ is not a reticulation, so there exists a node $u'\in N'\setminus Y'$ with children $x',y'\in Y'$. Let $x,y$ be the nodes in $Y$ corresponding to $x'$ and $y'$. The nodes $x$ and $y$ do not have a common parent not in $Y$, Let $z_x$ and $z_y$ be parents of $x$ and $y$ not contained in $Y$. Neither $(z_x,x)$ nor $(z_y,y)$ is movable, so $z_x$ forms a triangle together with its parent $c_x$ and its other child $d_x$, and similarly $(c_y,z_y,d_y)$ also forms a triangle.}
We distinguish three subcases:

%\todoCeline{Please use the same notation as Fig 17 of the main paper. {\bf RJ:} I think I do now, but the first part does not completely coincide with any previous case: they are subcases of 3biii}
\begin{enumerate}
\item\label{case:reorientate} {\bf $\bm{d_y\not\in Y}$.} Inverting the orientation of $(z_y,d_y)$ gives a valid orientation with the same underlying unrooted network, and conserves the downward-closedness of $Y$ and the isomorphism between $Y$ and $Y'$ (Figure~\ref{fig:ReorientateTriangle}). The resulting rooted network has a reticulation node directly above $Y$, so we can add a node to $Y$ and $Y'$ with at most one move.

\begin{figure}[h!]
\begin{center}
\includegraphics[scale=0.5]{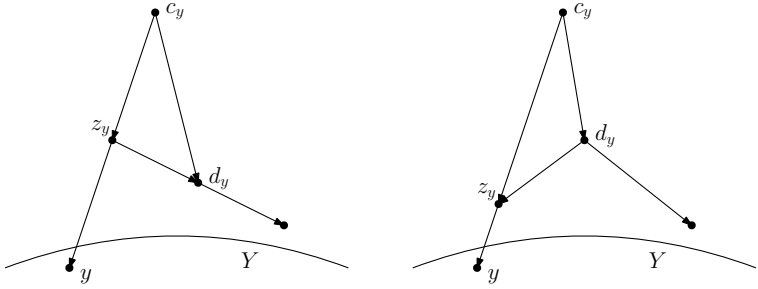}
\end{center}
\caption{The re-orientation of the bottom edge of a triangle, all of whose nodes are above $Y$.}
\label{fig:ReorientateTriangle}
\end{figure}

\item {\bf $\bm{d_x\not\in Y}$.} \celine{This case is handled symmetrically to Case~\ref{case:reorientate} by inverting the orientation of $(z_x,d_x)$}. %By symmetric arguments to Case~\ref{case:reorientate} we can add one node after re-orientating an edge and doing one move. }

% Now suppose that all nodes directly above $Y$ and $Y'$ are tree nodes, and that none of them are on the side of a triangle, of which the bottom node is not in $Y$ or $Y'$. By assumption, none of the incoming edges of $x'$ and $y'$ are movable, and then, by assumption, all are on the side of a triangle of which the bottom node is in $Y'$. We pick one of these bottom nodes $z'$ and note that both incoming edges are movable. Let $z\in Y$ be the corresponding node of $z'$ which has unmovable incoming edges by assumption. This means its parent edges are on the side of triangles, of which the bottom node is in $Y$. Pick such a bottom node $w$ and let $w'$ be the corresponding node in $Y'$. Now we can move an incoming edge of $z'$ to an incoming edge of $w'$ to create a tree node with children $w'$ and $z'$, which has existing corresponding tree node with children $w$ and $z$ in the other network. We can therefore add one node to $Y$ and $Y'$ using one move.

\item {\bf $\bm{d_x,d_y\in Y}$.} Note that one of $(u',x')$ and $(u',y')$ is movable \leo{(Observation~\ref{rem:MovableTree})}. 
\begin{enumerate}
\item {\bf $\bm{(u',x')}$ is movable.\label{case3}} 
%Let $x$ be the node in $Y$ corresponding to $x'$, and note that its parent $z_x$ is on the side of a triangle, with bottom node $d_x\in Y$. 
Tail moving $(u',x')$ to an incoming edge of $d_x'$, the node in $Y'$ corresponding to $d_x$, we create a lowest tree node $z_x'$ with children $x'$ and $d_x'$ in $Y'$ (Figure~\ref{fig:UnrootedOneMove}). We can add $z_x$ (with children $x$ and $d_x$) and $z_x'$ (with children $x'$ and $d_x'$) to $Y$ and $Y'$ at the cost of one move.
\item {\bf $\bm{(u',y')}$ is movable.} This case is handled symmetrically to the previous case  by  \leo{interchanging the roles of} $x$ and $y$.
\end{enumerate}

\end{enumerate}
\begin{figure}[h!]
\begin{center}
\includegraphics[scale=0.5]{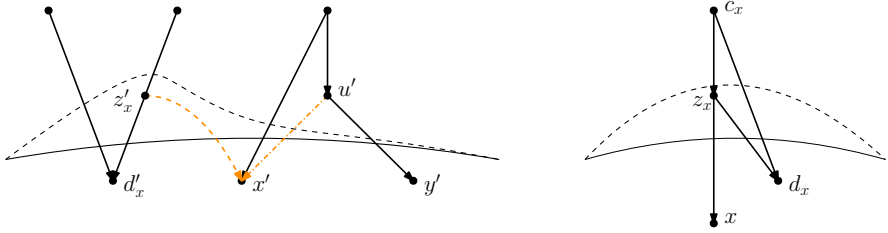}
\end{center}
\caption{The \celine{tail move used in the proof of Case \ref{case3} in the proof of Theorem \ref{thm:proofNNIimproved}}%one move used to add a node in the case that all nodes directly above $Y$ are tree nodes, and  none of these is on the side of a triangle whose bottom node is not in $Y$.
}
\label{fig:UnrootedOneMove}
\end{figure}

This proves that we can always add a node to $Y$ and $Y'$ using at most one SPR move. Hence the number of SPR moves needed to transform $N_r$ into $N_r'$ is at most \markj{$|V| = n+2k$}. Together with the moves needed to get to $N_r$ and $N_r'$ from $N$ and $N'$, we get an upper bound of \leo{$n+\frac{8}{3}k$} for the number of SPR moves needed to go from $N$ to $N'$.\qed
\end{proof}

\section{Discussion}

In a rooted phylogenetic tree, each rSPR move is a tail move, while in a network one can also perform head moves. Hence, it is natural to define rSPR moves in networks as the union of head and tail moves~\citep{Gambette-vI-Rearrangement}. However, we have shown that to connect the tiers of phylogenetic network space, tail moves are sufficient (except for one trivial case with only two taxa). \leo{In fact, tail moves, by themselves, can also be seen as a natural generalization of rSPR moves on trees to networks.}

%\todoRemie{I'd say both Tail moves and rSPR moves are natural generalisations of moves: Tail moves are the network version of rSPR on trees; rSPR moves on networks are the rooted version of SPR moves on networks.}

Nevertheless, there can be several reasons to consider head moves in addition to tail moves when exploring the tiers of network space. First of all, the number of tail moves that one needs to transform one network into another may be bigger than the number of rSPR moves. The difference is not very big though; we have shown that at most four tail moves are sufficient to mimic one head move (assuming~$|X|\geq 3$). Another reason to use head moves is \leo{when} a search using \leo{only} tail moves gets stuck in a local optimum. A head move can then be used to escape from this local optimum. Notice that a head move can basically move a reticulation to a completely different location in the network. Hence, from a biological point of view, head moves \leo{can} change a network more drastically than tail moves.

On rooted phylogenetic trees, rNNI moves can be seen as distance-1 rSPR moves. Therefore, it is natural to define rNNI moves in the same way on networks \citep[see][]{Gambette-vI-Rearrangement}. Hence, rNNI moves on networks are distance-1 head moves and distance-1 tail moves. Also for this case we have shown that the distance-1 head moves \leo{may be omitted while preserving} connectivity (again excluding one two-taxa case).

To determine which moves explore network space most efficiently, practical experiments are necessary. It seems reasonable to believe that a combination of \leo{different types of} moves will work best. Distance-1 head and tail moves moves are local moves that change only a small part of the network and the number of possible moves is relatively small. Therefore, such moves are suitable for finding a local optimum. General tail and head moves are not local and can therefore be useful to move away from local optima in order to search for a global optimum. Implementing and analysing such local-search methods is an important topic for future research.

\begin{figure}[h!]
\begin{center}
\includegraphics[scale=0.55]{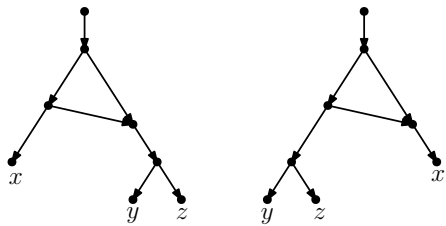}
\end{center}
\caption{Two networks for which subtree reduction does not preserve tail-move distance.}
\label{fig:SubtreeReductionFail}
\end{figure}

Another interesting direction is to develop algorithms for computing the tail-move distance between networks. Although we have shown that this problem is NP-hard, we do not know, for example, whether it is fixed-parameter tractable, with the tail move distance as parameter. We note that several common reduction rules, such as subtree and cluster reduction are not always safe for this problem (see \citep{bordewich2005computational,bordewich2017fixed} for a definition of these reductions and the chain reduction mentioned below). As an extreme example, consider the networks in Figure~\ref{fig:SubtreeReductionFail}, which have tail-move distance~3. After reducing the subtrees on $\{y,z\}$ to a single leaf, the tail move distance becomes infinite. It could, however, be safe to reduce subtrees to size~2. Similar arguments hold for cluster reduction. Another interesting question is whether reducing chains to length~2 (or any other constant length) is safe.%\todo{The networks in Figure~\ref{fig:SubtreeReductionFail} need a root edge!}

Finally, a major open problem concerns the rNNI diameter as well as the distance-1 tail move diameter: Are these diameters quadratic in the number of leaves, or $O(n\log(n))$? For \leo{tail, rSPR and SPR} moves, we know that the diameter is linear in the number of leaves, but it would be interesting to \leo{find tight bounds (see Table~\ref{tab:diameters})}.\\

%\todoCeline{Discuss locally of moves distance-1 tail, tail, rSPR, NNI}
%Open:
% algorithm computing exact distance
% computation FPT?
% Sensible kernelisation
% subtree, cluster and chain reduction to size-2 are safe?
% Diameter bound improvements: 
%  NNI/Tail_1 lower/upper bound (order nlog(n) or n^2?)
%  SPR/Tail constants
% 1 leaf more optimal move (d)?

\noindent\small{{\bf Funding} Remie Janssen and Mark Jones were supported by Vidi grant 639.072.602 from the Netherlands Organization for Scientific Research (NWO). Leo van Iersel was partly supported by NWO, including Vidi grant 639.072.602, and partly by the 4TU Applied Mathematics Institute. P\'eter L. Erd\H{o}s was partly supported by the Hungarian National Research, Development and Innovation Office NKFIH, under the grants K 116769 and SNN 116095. Celine Scornavacca was partially supported  by the French Agence Nationale de la Recherche Investissements d'Avenir/ Bioinformatique (ANR-10-BINF-01-02, Ancestrome).}

\bibliographystyle{chicago}

\bibliography{bibliography}

\begin{thebibliography}{}

\bibitem[\protect\citeauthoryear{Abbott, Albach, Ansell, Arntzen, Baird,
  Bierne, Boughman, Brelsford, Buerkle, Buggs, et~al.}{Abbott
  et~al.}{2013}]{Abbott}
Abbott, R., D.~Albach, S.~Ansell, J.~W. Arntzen, S.~J. Baird, N.~Bierne,
  J.~Boughman, A.~Brelsford, C.~A. Buerkle, R.~Buggs, et~al. (2013).
\newblock Hybridization and speciation.
\newblock {\em Journal of Evolutionary Biology\/}~{\em 26\/}(2), 229--246.

\bibitem[\protect\citeauthoryear{Atkins and McDiarmid}{Atkins and
  McDiarmid}{2015}]{atkins2015extremal}
Atkins, R. and C.~McDiarmid (2015).
\newblock Extremal distances for subtree transfer operations in binary trees.
\newblock {\em arXiv preprint arXiv:1509.00669\/}.

\bibitem[\protect\citeauthoryear{Bordewich, Linz, and Semple}{Bordewich
  et~al.}{2017}]{Bordewich-Lost}
Bordewich, M., S.~Linz, and C.~Semple (2017).
\newblock Lost in space? {G}eneralising subtree prune and regraft to spaces of
  phylogenetic networks.
\newblock {\em Journal of Theoretical Biology\/}~{\em 423}, 1--12.

\bibitem[\protect\citeauthoryear{Bordewich, Scornavacca, Tokac, and
  Weller}{Bordewich et~al.}{2017}]{bordewich2017fixed}
Bordewich, M., C.~Scornavacca, N.~Tokac, and M.~Weller (2017).
\newblock On the fixed parameter tractability of agreement-based phylogenetic
  distances.
\newblock {\em Journal of mathematical biology\/}~{\em 74\/}(1-2), 239--257.

\bibitem[\protect\citeauthoryear{Bordewich and Semple}{Bordewich and
  Semple}{2005}]{bordewich2005computational}
Bordewich, M. and C.~Semple (2005).
\newblock On the computational complexity of the rooted subtree prune and
  regraft distance.
\newblock {\em Annals of combinatorics\/}~{\em 8\/}(4), 409--423.

\bibitem[\protect\citeauthoryear{Ding, Gr{\"u}newald, and Humphries}{Ding
  et~al.}{2011}]{ding2011agreement}
Ding, Y., S.~Gr{\"u}newald, and P.~J. Humphries (2011).
\newblock On agreement forests.
\newblock {\em Journal of Combinatorial Theory, Series A\/}~{\em 118\/}(7),
  2059--2065.

\bibitem[\protect\citeauthoryear{Felsenstein}{Felsenstein}{2004}]{felsenstein2004inferring}
Felsenstein, J. (2004).
\newblock {\em Inferring phylogenies}, Volume~2.
\newblock Sinauer associates Sunderland, MA.

\bibitem[\protect\citeauthoryear{Francis, Huber, Moulton, and Wu}{Francis
  et~al.}{2017}]{francis2017bounds}
Francis, A., K.~Huber, V.~Moulton, and T.~Wu (2017).
\newblock Bounds for phylogenetic network space metrics.
\newblock {\em arXiv preprint arXiv:1702.05609\/}.

\bibitem[\protect\citeauthoryear{Gambette, van Iersel, Jones, Lafond, Pardi,
  and Scornavacca}{Gambette et~al.}{2017}]{Gambette-vI-Rearrangement}
Gambette, P., L.~van Iersel, M.~Jones, M.~Lafond, F.~Pardi, and C.~Scornavacca
  (2017).
\newblock Rearrangement moves on rooted phylogenetic networks.
\newblock {\em Accepted for publication in PLOS Computational Biology\/}.

\bibitem[\protect\citeauthoryear{Huber, Moulton, and Wu}{Huber
  et~al.}{2016}]{huber2016transforming}
Huber, K.~T., V.~Moulton, and T.~Wu (2016).
\newblock Transforming phylogenetic networks: Moving beyond tree space.
\newblock {\em Journal of Theoretical Biology\/}~{\em 404}, 30--39.

\bibitem[\protect\citeauthoryear{Huson, Rupp, and Scornavacca}{Huson
  et~al.}{2010}]{Huson}
Huson, D.~H., R.~Rupp, and C.~Scornavacca (2010).
\newblock {\em Phylogenetic networks: concepts, algorithms and applications}.
\newblock Cambridge University Press.

\bibitem[\protect\citeauthoryear{Lempel, Even, and Cederbaum}{Lempel
  et~al.}{1967}]{lempel1967algorithm}
Lempel, A., S.~Even, and I.~Cederbaum (1967).
\newblock An algorithm for planarity testing of graphs.
\newblock In {\em Theory of graphs: International symposium}, Volume~67, pp.\
  215--232. Gordon and Breach, New York.

\bibitem[\protect\citeauthoryear{Li, Tromp, and Zhang}{Li
  et~al.}{1996}]{li1996nearest}
Li, M., J.~Tromp, and L.~Zhang (1996).
\newblock On the nearest neighbour interchange distance between evolutionary
  trees.
\newblock {\em Journal of Theoretical Biology\/}~{\em 182\/}(4), 463--467.

\bibitem[\protect\citeauthoryear{Morrison}{Morrison}{2011}]{Morrison}
Morrison, D. (2011).
\newblock {\em Introduction to Phylogenetic Networks}.
\newblock RJR Productions.

\bibitem[\protect\citeauthoryear{Than, Ruths, and Nakhleh}{Than
  et~al.}{2008}]{Than}
Than, C., D.~Ruths, and L.~Nakhleh (2008).
\newblock Phylonet: a software package for analyzing and reconstructing
  reticulate evolutionary relationships.
\newblock {\em BMC bioinformatics\/}~{\em 9\/}(1), 322.

\bibitem[\protect\citeauthoryear{Vuilleumier and Bonhoeffer}{Vuilleumier and
  Bonhoeffer}{2015}]{Vuilleumier}
Vuilleumier, S. and S.~Bonhoeffer (2015).
\newblock Contribution of recombination to the evolutionary history of hiv.
\newblock {\em Current Opinion in HIV and AIDS\/}~{\em 10\/}(2), 84--89.

\bibitem[\protect\citeauthoryear{Yu, Dong, Liu, and Nakhleh}{Yu
  et~al.}{2014}]{yu2014maximum}
Yu, Y., J.~Dong, K.~J. Liu, and L.~Nakhleh (2014).
\newblock Maximum likelihood inference of reticulate evolutionary histories.
\newblock {\em Proceedings of the National Academy of Sciences\/}~{\em
  111\/}(46), 16448--16453.

\bibitem[\protect\citeauthoryear{Zhang, Ogilvie, Drummond, and Stadler}{Zhang
  et~al.}{2017}]{Zhang124982}
Zhang, C., H.~A. Ogilvie, A.~J. Drummond, and T.~Stadler (2017).
\newblock Bayesian inference of species networks from multilocus sequence data.
\newblock {\em bioRxiv\/}.

\bibitem[\protect\citeauthoryear{Zhaxybayeva and Doolittle}{Zhaxybayeva and
  Doolittle}{2011}]{Zhaxybayeva}
Zhaxybayeva, O. and W.~F. Doolittle (2011).
\newblock Lateral gene transfer.
\newblock {\em Current Biology\/}~{\em 21\/}(7), R242--R246.

\end{thebibliography}

\end{document}